\documentclass[12pt]{amsart}
\usepackage{amsmath}
\usepackage{amscd}
\usepackage{amssymb}
\usepackage{amsfonts}
\usepackage{amsthm}
\usepackage{bbm}
\usepackage{cancel}
\usepackage{color}
\usepackage{eucal}
\usepackage{enumerate,yfonts}
\usepackage{enumitem}
\usepackage{pdfsync}
\usepackage[all,cmtip]{xy}
\usepackage{graphicx}
\usepackage{graphics}
\usepackage{hyperref}
\usepackage{latexsym}
\usepackage{mathrsfs}
\usepackage{placeins}
\usepackage{pstricks}
\usepackage{bookmark}
\usepackage{mathtools}
\usepackage{stmaryrd}
\usepackage{url}
\usepackage{mathrsfs}
\usepackage{amsmath,calligra,mathrsfs}
\DeclareMathOperator{\Hom}{\mathscr{H}\text{\kern -3pt {\calligra\large om}}\,}

\newtheorem{thm}{Theorem}[section]
\newtheorem{theorem}[thm]{Theorem}

\newtheorem{lemma}[thm]{Lemma}

\newtheorem{proposition}[thm]{Proposition}

\newtheorem{thm-dfn}[thm]{Theorem-Definition}

\newtheorem{remark}[thm]{Remark}

\numberwithin{equation}{section}

\newenvironment{rouge}
{\relax\color{red}}
{\hspace*{.3ex}\relax}

\newenvironment{bluet}
{\relax\color{blue}}
{\hspace*{.3ex}\relax}

\newcommand{\br}{\begin{rouge}}
\newcommand{\er}{\end{rouge}}

\newcommand{\bb}{\begin{bluet}}
\newcommand{\eb}{\end{bluet}}

\newcommand{\nc}{\newcommand}

\newcommand{\cI}{{\mathcal I}}

\newcommand{\cD}{{\mathcal D}}
\newcommand{\cO}{{\mathcal O}}

\newcommand{\cF}{{\mathcal F}}
\newcommand{\cN}{{\mathcal N}}
\newcommand{\cH}{{\mathcal H}}

\newcommand{\cL}{{\mathcal L}}

\newcommand{\cZ}{{\mathcal Z}}

\newcommand{\cM}{{\mathcal M}}
\newcommand{\cG}{{\mathcal G}}

\newcommand{\bC}{{\mathbb C}}
\newcommand{\bZ}{{\mathbb Z}}

\newcommand{\bQ}{{\mathbb Q}}
\newcommand{\bD}{{\mathbb D}}
\newcommand{\bN}{{\mathbb N}}

\newcommand{\bfV}{\mathbf V}

\newcommand{\rr}{\mathrm{r}}

\newcommand{\rrs}{\mathrm{rs}}

\newcommand{\La}{{\mathfrak{a}}}

\newcommand{\Lg}{{\mathfrak{g}}}

\newcommand{\Ll}{{\mathfrak{l}}}

\newcommand{\Ls}{{\mathfrak{s}}}
\newcommand{\Lu}{{\mathfrak{u}}}

\newcommand{\fg}{{\mathfrak g}}

\newcommand{\fn}{{\mathfrak{n}}}

\newcommand{\fl}{{\mathfrak{l}}}

\newcommand{\fm}{{\mathfrak{m}}}
\newcommand{\fk}{{\mathfrak{K}}}

\newcommand{\fa}{{\mathfrak{a}}}

\nc{\ot}{\otimes}
\nc{\on}{\operatorname}

\nc{\oh}{{\operatorname{H}}}
\nc{\gr}{{\operatorname{gr}}}
\nc{\Gr}{{\operatorname{Gr}}}
\nc{\rk}{{\operatorname{rank}}}
\nc{\codim}{{\operatorname{codim}}}
\nc{\img}{{\operatorname{Im}}}
\nc{\IC}{{\operatorname{IC}}}
\nc{\CC}{{\operatorname{CC}}}

\nc{\bI}{{\mathbf 1}}

\nc{\lp}{{\left(}}
\nc{\rp}{{\right)}}

\newcommand{\beqn}{\begin{equation*}}
\newcommand{\eeqn}{\end{equation*}}

\newcommand{\beq}{\begin{equation}}
\newcommand{\eeq}{\end{equation}}

\newcommand{\bega}{\begin{gathered}}
\newcommand{\eega}{\end{gathered}}

\newcommand{\bern}{\begin{eqnarray*}}
\newcommand{\eern}{\end{eqnarray*}}
\newcommand{\ber}{\begin{eqnarray}}
\newcommand{\eer}{\end{eqnarray}}

\newcommand{\inv}{{\mathbin{/\mkern-4mu/}}}

\newcommand{\barbasepta}{\bar{a}_0}

\newcommand{\Cartan}{\fa}
\newcommand{\Sym}{\on{Sym}}
\newcommand{\Mod}{\on{Mod}}
\newcommand{\Char}{\mathscr{C}\text{\kern -3pt {\calligra har}}}
\nc{\del}{{\partial }}
\nc{\reg}{{\mathrm{sr}}}

\nc{\Ts}{{T_{\mathbf{s}}}}
\nc{\Tds}{{\check T_{\mathbf{s}}}}
\nc{\Lie}{{\operatorname{Lie}}}

\nc{\ke}{k+e}

\begin{document}

\title[Invariant systems and character sheaves for graded Lie algebras]{Invariant systems and character sheaves for graded Lie algebras}

\author{Kari Vilonen}\address{School of Mathematics and Statistics, University of Melbourne, VIC 3010, Australia, also Department of Mathematics and Statistics, University of Helsinki, Helsinki, Finland}
\email{kari.vilonen@unimelb.edu.au, kari.vilonen@helsinki.fi}
\thanks{KV was supported in part by the ARC grants DP150103525, DP180101445,  FL200100141 and the Academy of Finland}

\author{Ting Xue}
\address{School of Mathematics and Statistics, University of Melbourne, VIC 3010, Australia, also Department of Mathematics and Statistics, University of Helsinki, Helsinki, Finland} 
\email{ting.xue@unimelb.edu.au}
\thanks{TX was supported in part by the ARC grant DP150103525.}

\begin{abstract}
In this paper we introduce a new ingredient, invariant systems of differential equations, to our study of character sheaves on graded Lie algebras. The character sheaves we construct in this paper, together with the ones constructed in~\cite{VX2,X2}, constitute all cuspidal character sheaves for Vinberg's type I classical graded Lie algebras. 
\end{abstract}


\maketitle

\tableofcontents

\section{Introduction}

In this paper we continue our study of character sheaves on graded Lie algebras initiated in~\cite{VX2}. The goal is to explicitly determine all cuspidal character sheaves, i.e., the irreducible perverse sheaves whose singular support is nilpotent and which do not arise by parabolic induction from smaller graded Lie algebras. In the ungraded case this problem was solved by Lusztig when he established the generalized Springer correspondence~\cite{Lu}. Combining the results of this paper with those of~\cite{LTVX} and~\cite{VX2} we accomplish this goal in majority of cases. We will make the word ``majority" more precise later in this introduction. 

Our main geometric tool is the nearby cycle construction initially due to Grinberg in~\cite{G2} and further studied in more generality suitable for our purposes in~\cite{GVX1,GVX2}. The main result of~\cite{LTVX} states that all cuspidal character sheaves can be obtained from this construction in a very specific manner. 

The discussion above reduces the question of explicitly describing all the cuspidal character sheaves to a detailed study of the nearby cycle construction. In~\cite{VX2} we studied the (GIT) stable gradings using~\cite{GVX2}. In this paper we treat the general case allowing unstable gradings. As the unstable case is substantially more complicated, we were not able to carry out the study with geometric tools only. Thus, we introduce a new technique, invariant systems of differential equations, generalising the classical Harish-Chandra system. In the ungraded case  Hotta and Kashiwara showed in~\cite{HK} that the Harish-Chandra system coincides with the Springer sheaf, which in turn coincides with the sheaf obtained from the nearby cycle construction. 
However, in the graded case this invariant system very rarely coincides with the nearby cycle sheaf. In this paper we show, in  \autoref{sec-rank1g}, that under very strong conditions in the case of rank one graded Lie algebras, the invariant system and the sheaf obtained by the nearby cycle  construction  are closely related.

The advantage of the invariant system is that it can readily be analysed by solving the corresponding system of differential equations. A crucial ingredient in allowing us to analyse the system stems from the fact that invariant systems are controlled by certain b-functions. Moreover, in the cases of interest to us the b-functions are known explicitly thanks to the work of Fumihiro Sato~\cite{Saf}. The use of the invariant systems lets us read off certain Hecke relations which are crucial in determining the cuspidal character sheaves.

In this paper we  exhibit all cuspidal character sheaves for graded Lie algebras in explicit terms, in the case of classical graded Lie algebras. Recall that in~\cite{Vin} Vinberg divides classical graded Lie algebras in three different types, called type I, II, and III. It can be deduced from the results in~\cite{LTVX}, as conjectured in~\cite{X2}, that in type II  there are no cuspidal character sheaves, and that  in type III cuspidal character sheaves have nilpotent support and are rare, just as in the ungraded situation~\cite{Lu} and in the $\bZ$-graded case~\cite{Lu2}. Similarly one can deduce from~\cite{LTVX}
that the character sheaves constructed in this paper in  \autoref{cuspidal sheaves}, combined with those in\cite{VX2,X2}, constitute all cuspidal character sheaves for classical  graded Lie algebras of type I. 

In what follows we explain the content of this paper in a bit more detail. Let $G$ be a complex semisimple algebraic group and $\Lg$ its Lie algebra. Let $\theta$ be an automorphism of $G$ of finite order $m$. It induces a $\bZ/m\bZ$-grading on the Lie algebra 
\beqn
\Lg = \bigoplus_{i\in\bZ/m\bZ} \fg_i
\eeqn
 and $K:=(G^\theta)^0$ acts on the constituents $\Lg_i$ by adjoint action.  We can identify $\fg_1^*$ with $\fg_{m-1}=\fg_{-1}$. 
Recall that $\Lg_1$ possesses a Cartan subspace $\fa\subset \fg_1$ which is a maximal abelian subspace consisting of semisimple elements. We have the adjoint quotient map~\cite{Vin}
\beqn
f: \Lg_1 \to \Lg_1 \inv K \cong \La \slash W_\La \,,
\eeqn
where $W_\fa=N_K(\fa)/Z_K(\fa)$, the little Weyl group, is a complex reflection group.  By the rank of the graded Lie algebra we mean the dimension of the Cartan subspace $\fa$, and by the {\em semisimple} rank we mean $\dim\fa-\dim\fa^K$. Note that we have a similar picture for $\Lg_{-1}$. Recall that, by definition a character sheaf on $\Lg_1$ is a simple $K$-equivariant perverse sheaf on $\Lg_1$ with nilpotent singular support. Put differently, they are Fourier transforms of simple $K$-equivariant perverse sheaves on the nilpotent cone $\cN_{-1}=\cN\cap\fg_{-1}$ in $\fg_{-1}$, where $\cN$ is the nilpotent cone of $\Lg$.

To make the  ``nearby cycle construction" we consider a general fibre $F$ of the adjoint quotient map $f^{-}:\Lg_{-1}\to\Lg_{-1}\inv K$. On the open $K$-orbit in $F$ we fix a clean $K$-equivariant local system $\cL_\chi$. We form the corresponding intersection cohomology sheaf $\IC(\cL_\chi)$, take its limit to the nilpotent cone $\cN_{-1}$ (the zero fibre of $f^-$),  and  then take the Fourier transform. We denote the limit by $P_\chi\in\on{Perv}_K(\cN_{-1})$ and its Fourier transform by $\widehat P_\chi\in\on{Perv}_K(\Lg_1)$. The irreducible constituents of $\widehat P_\chi$ are character sheaves by definition. It turns out that under our hypotheses the resulting sheaf $\widehat P_\chi$ is an intersection cohomology  sheaf  whose support is the entire $\fg_1$. We proceed as in~\cite{GVX2} to analyse the resulting local system on the generic locus of $\fg_1$ associated to $\widehat P_\chi$.  The local system can be expressed in terms of a Hecke algebra of a complex reflection group with certain parameters. The parameters give the Hecke relations and determining these Hecke relations is a key to understanding the local system.

In  \autoref{nearby} we show how the determination of the Hecke relations is reduced to the case of semisimple rank one. The determination of the Hecke relations in semisimple rank one stable cases was done topologically in~\cite{VX2}. This is feasible because the geometry in those cases is rather simple. The fundamental invariant $f$, i.e. the map $f: \Lg_1 \to  \La \slash W_\La \cong\bC$, is given by a monomial so that the null-cone is a (non-reduced) normal crossings divisor. In the unstable case the invariant $f$ is much more complicated and we do not know how to perform the calculation topologically. To perform this calculation we introduce an invariant system $\cM_\chi$. For the precise definition of $\cM_\chi$ see~\eqref{rk1 M}. Then in \autoref{sec-rank1g} we show that under very strong conditions on the zeros of certain b-functions the $\widehat P_\chi$ and the invariant system $\cM_\chi$ coincide generically and this allows us to use $\cM_\chi$ to obtain the Hecke relations from the b-functions. It is a small miracle that in the cases where  cuspidal character sheaves exist the strong conditions on b-functions hold. Furthermore, explicit formulas for these b-functions can be deduced from calculations of b-functions of Fumihiro Sato~\cite{Saf} based on the study of zeta functions associated to prehomogenous vector spaces. The fact that we can explicitly obtain the b-functions we need from~\cite{Saf} was pointed out to us by Akihito Wachi.

In the ungraded case Lusztig~\cite{Lu} shows that semisimple Lie algebras afford cuspidal character sheaves only under very restrictive conditions. The situation is similar in the case of graded Lie algebras, i.e., only under some conditions on the gradings can there be cuspidal character sheaves. In  \autoref{cuspidal sheaves}  we write down these (type I) classical graded Lie algebras in a very concrete form as representations of certain quivers with prescribed dimension vectors. By~\cite{LTVX} these are precisely the gradings for which cuspidal character sheaves exist.  In Theorem~\ref{main theorem} we then write down completely explicitly the cuspidal character sheaves for these gradings.

The paper is organised as follows.  In  \autoref{invariant} we give a general discussion of invariant systems for polar representations and also explain the nearby cycle construction in this situation. In  \autoref{zero} we discuss the case of rank zero, i.e., the case of prehomogenous vector spaces. In  \autoref{sec-rank1g} we analyse invariant systems of polar representations of rank one. This case is closely related to the theory of prehomogenous vector spaces. We show that under strong conditions on zeros of the b-functions, the invariant system is closely related to the $\cD$-module arising from the nearby cycle construction. In  \autoref{graded} we recall basic facts about graded Lie algebras and the associated equivariant geometry. In  \autoref{nearby} we explain how to modify the nearby cycle construction of~\cite{GVX2} so that it works in our situation to allow reduction to semisimple rank one. In \autoref{Exact} we spell out the translation of the results on invariant systems to the graded Lie algebras setting we consider. Finally, in  \autoref{cuspidal sheaves}, we apply the considerations in previous sections to explicitly construct cuspidal character sheaves for (type I) classical graded  Lie algebras. 

{\bf Acknowledgement.} We thank Misha Grinberg and Toby Stafford for helpful conversations. Special thanks go to Akihito Wachi who answered our questions about b-functions and directed us to the reference~\cite{Saf}. TX thanks Cheng-Chiang Tsai and Zhiwei Yun for many helpful discussions.

\section{Invariant systems and the nearby cycle construction}
\label{invariant}

In this section we consider two constructions for character sheaves in the general setting of polar representations of~\cite{DK}. The first one is invariant systems of differential equations. They can be considered as a straightforward generalization of the classical Harish-Chandra system. The second construction is via nearby cycles followed by the Fourier transform. It is a generalization of the classical Springer construction which was already considered in~\cite{G2} in less generality. In the classical case of the adjoint action of a complex reductive group on its Lie algebra these two constructions coincide as was shown by~\cite{HK}. However, for general polar representations this is no longer the case and the two constructions are quite different. Under very strong conditions  the two constructions (for rank one polar representations) are closely related as we will explain later in~\autoref{ssec-comparison}.

Throughout this paper we work with complex algebraic varieties. When the variety is smooth we work with algebraic $\cD$-modules which for us will always be holonomic and with regular singularities. We identify such $\cD$-modules with perverse sheaves via the de Rham functor. Consequently our  perverse sheaves are  in the setting of classical topology and with complex coefficients. We will shift between these two points of view freely. 

Although in this section we work in the general context of polar representations,  the specific application we have in mind is the context of graded Lie algebras. One can substitute the graded Lie algebras setup of  \autoref{graded} for polar representations while reading this and the two subsequent sections.  The context of polar representations helps isolate the properties that are essential for our constructions although at times we impose strong conditions on them.

\subsection{Visible polar representations}Let $V$ be a polar representation of a complex  reductive algebraic group $G$ as introduced in \cite{DK}. 
Let $G_v$ denote the stabiliser subgroup in $G$ for $v\in V$. We say that $v\in V$ is semisimple if the orbit $G\cdot v$ is closed. We call $v\in V$ {\em regular} if $\on{dim}G_v\leq\on{dim}G_w$ for all $w\in V$. We call a semisimple element $v\in V$ {\em regular semisimple} if $\on{dim}G_v\leq\on{dim}G_w$ for all {\em semisimple} $w\in V$. Note that a regular semisimple element is not necessarily regular, but is only regular among semisimple elements. 

A key property of a polar representation $G|V$ is that there exists a Cartan subspace $\fa\subset V$ consisting of semisimple elements such that $\dim\fa=\dim V\inv G$. Let $\Cartan \subset V$ be a Cartan subspace (all Cartan subspaces are conjugate under $G$) and  let $W = N_G (\Cartan) / Z_G (\Cartan)$ be the associated Weyl group. Then $V\inv G \cong \Cartan/W$ and, when $G$ is connected, $W$ is a  complex reflection group acting on $\Cartan$   (see~\cite[Theorems 2.9, 2.10]{DK}). The rank of $(G,V)$ is defined to be the dimension of $\fa$. 

We write $$Q = \Cartan / W = V \inv G.$$ Let $f : V \to Q$ be the quotient map. The representation $G | V$ is called {\it visible} if the null-cone $\mathfrak{N}=f^{-1} (0)$ consists of finitely many $G$-orbits. From now on we assume visibility. This condition implies that any fibre of the map $f$ consists of finitely many $G$-orbits. Let $Q^{\rrs}=\Cartan^\rrs/W$, where $\Cartan^\rrs$ is the set of regular semisimple elements in $\fa$. We call an element $v\in V$ strongly regular if $v$ is regular and $f(v)\in Q^{\rrs}$.  We write $V^\reg$ for the open subset of strongly regular elements in $V$. 

Consider a generic fibre $F=f^{-1}(q)$, for $q\in Q^{\rrs}$. We analyse the fibre $F$ for later use.  Let $Z_G(\Cartan)^0$ denote the identity component of the centraliser of $\Cartan$ in $G$. We have the following $Z_G(\Cartan)^0$-invariant decomposition
\begin{equation}
\label{decomp}
V\ = \ U \oplus \Cartan \oplus \fg\cdot \Cartan
\end{equation}
where $U$ is a $Z_G(\Cartan)^0$-invariant complement of $\Cartan \oplus \fg\cdot \Cartan$. The $U$ consists of finitely many orbits under $Z_G(\Cartan)^0$-action, i.e., it is a prehomogenous vector space for $Z_G(\Cartan)^0$. We write $U^\mathrm{r}$ for the open $Z_G(\Cartan)^0$-orbit in $U$. It consists of regular elements for the polar representation $(Z_G(\Cartan)^0,U)$.  
Then we have
\beq
\label{reg decomp}
V^{\reg}\ = \ G\cdot(U^\mathrm{r} \times\Cartan^\mathrm{rs})\,. 
\eeq
We also see that
\begin{equation}\label{intersection}
F\cap \Cartan = W\cdot v\,,
\end{equation}
where $v$ is any element in $F\cap \Cartan \subset  \Cartan^\rrs$.

\begin{remark}When $U=0$ in~\eqref{decomp}, the polar representation is called stable in~\cite{DK}. Note that this notion is different from GIT stability. For stable polar representations, we have $V^{\reg}=G\cdot\Cartan^\mathrm{rs}$.
\end{remark}

 Let us write $f_1, \dots , f_r$ for the fundamental invariants so that $\Sym(V^*)^G$ is a polynomial algebra in the  $f_1, \dots , f_r$. Let us dually consider $\Sym(V)^G$ and write $D_1, \dots , D_r$ for the generators which are the duals of the  $f_1, \dots , f_r$. We view $D_1, \dots , D_r$ as (invariant) constant coefficient differential operators on $V$ and as such as elements in the Weyl algebra $\cD_V$ on $V$. 

\subsection{Invariant systems}Let $\Mod^G_{hr}(\cD_V)$ denote the category of $G$-equivariant regular  holonomic $\cD_V$-modules on $V$. We consider the subcategory $\Char_G(V)$ of $\Mod^G_{hr}(\cD_V)$ consisting of modules whose singular support is nilpotent, i.e., it is contained in $V\times \mathfrak{N}^*$, where $ \mathfrak{N}^*$ is the null-cone for the dual representation $V^*$. As in our previous work, we call the irreducible objects in $\Char_G(V)$  character sheaves on $V$. The character sheaves are thus Fourier transforms of $G$-equivariant $\cD$-modules on $\mathfrak{N}^*$. As, by visibility, $G$ has finitely many orbits on $\mathfrak{N}^*$, we see that $\cD$-modules whose singular support is nilpotent are automatically regular holonomic. As the objects in $\Mod^G_{hr}(\cD_{\mathfrak{N}^*})$ are monodromic under the $\bC^*$-action, the same is true for character sheaves.

Consider a $G$-equivariant cyclic module $\cM = \cD_V u=\cD_V/\cI$. Such a module lies in $\Char_G(V)$  if and only if the symbol ideal $\sigma(\cI)$ satisfies:
\beq\label{cond-char}
\sigma(\cI) \supset (\Sym(V)^G_+)^k \qquad \text{for some $k$.}
\eeq 
Differentiating the  $G$-action $\rho:G \to GL(V)$ gives rise to a representation $d\rho:\fg \to \on{End}(V)$.  Consider the following invariant system on $V$ associated to  a homomorphism $\phi:\fg\to\bC$:
\beq
\label{HC analogue}
\cM_\phi \ = \ \cD_V / \langle D_1,\dots, D_r , \{d\rho(\xi)-\phi(\xi) \mid \xi\in \fg \}\rangle\,.
\eeq
The system $\cM_\phi$ is regular holonomic.  
For general  homomorphisms $\phi$ the module $\cM_\phi$ is not necessarily $G$-equivariant but  is only $\phi$-monodromic. Considering such systems goes back at least to~\cite{Ho}. If $\phi$ is the differential of a character $\tilde\phi:G \to \bC^*$ then the resulting module $\cM_\phi$ is $G$-equivariant. It is not difficult to determine when two such characters give rise to the same equivariant module. We will analyse these systems in more detail in the cases of rank zero and rank one polar representations  in the next two sections. Such systems are obvious generalisations of Harish-Chandra systems. See~\cite{HK}, for example, where the relationship of the Harish-Chandra systems to Springer theory is  explained. 

Let now $\phi=d\tilde \phi$ where $\tilde\phi:G \to \bC^*$ is a character. Note that ${\cM_\phi}|_{V^\reg}$ is a vector bundle with a flat connection, of rank  $|W|=\prod\on{deg}D_i$.  The intermediate extensions of the irreducible constituents of $\cM|_{V^\reg}$ are character sheaves by definition (see~\eqref{cond-char}). 

\begin{remark}\normalfont
In~\cite{BNS} the authors also consider invariant systems of the form~\eqref{HC analogue} in the context of polar representations. They impose strong conditions on the polar representations and under those conditions obtain results similar to Harish-Chandra's regularity theorem. The invariant systems of interest to us are, in general, far from satisfying the conditions in~\cite{BNS} as we do not assume stability and do not assume that the image of the radial map is simple. We impose different strong conditions in the case of rank one polar representations so that the  nearby cycle construction and the invariant system are closely related. We will show that under those conditions the resulting sheaf of the nearby cycle construction can be described using the corresponding invariant system. 
\end{remark}

\subsection{Nearby cycles and their Fourier transform}We do not know how one can determine the character sheaves using invariant systems directly. Thus, we now turn to the nearby cycle construction which can be considered as a reasonable analogue of the Springer theory. We expect that this construction (and its slight generalisation) is sufficient to determine all character sheaves in the context of graded Lie algebras. 

We begin with a mild generalisation of the constructions in~\cite{G2,GVX1,GVX2}.  
As in~\cite[\S 2.5]{GVX2} we consider a generic fibre $F=f^{-1}(q)$, for $q\in Q^{\rrs}$ and restrict the family $f: V \to Q$ as follows:
\beq\label{eqn-family-Z}
f_q:\cZ= \{ (x, c) \in V \times \bC \; | \; f (x) = c \, q \} \to \bC\,,\,(x,c)\mapsto c, 
\eeq
which is by construction $G$-equivariant. We have a $G$-equivariant identification $f_q^{-1}(\bC^*) \cong F \times \bC^*$. Thus, given any $\cF \in \Mod^G_{hr}(\cD_F)$ we can regard it as a $\cD$-module on $f_q^{-1}(\bC^*)$ and then form the nearby cycle module $P_\cF=\psi_{f_q}(\cF)\in\Mod^G_{hr}(\cD_V)$ supported on the null-cone $f^{-1}(0)$.  
 Consider the Fourier transform $\widehat {P_\cF}$. As we are in the $\bC^*$-monodromic situation $\widehat {P_\cF}$ is also regular holonomic.  Hence,  by construction, $\widehat {P_\cF}\in \Char_G(V^*)$ and so its irreducible constituents are character sheaves on $V^*$.

This construction a priori depends on the point $q\in Q^{\rrs}$.   To analyse the dependence of $\psi_{f_q}(\cF)$ on $q$ we pass to the topological side. The key point is that the $G$-equivariance and visibility imply that the map $f$ satisfies Thom's condition $A_f$. This allows us to equip $\psi_{f_q}(\cF)$ with an explicit large automorphism group and, in particular, conclude that the isomorphism class of $\psi_{f_q}(\cF)$  does not depend on $q$.

 We will now explain how to determine  $\widehat {P_\cF}|_{(V^*)^\reg}$ topologically.  
In the following discussion we fix $v\in F\cap \Cartan\subset \Cartan^\rrs$ (see~\eqref{intersection}), but it is easy to see that this choice is not essential.  The fiber $F$ consists of finitely many $G$-orbits. The open orbit $F^\reg$ consists of strongly regular elements in $F$,  and the closed orbit $F^\mathrm{ss}$ consists of regular semisimple elements in $F$.  

We have
\begin{equation}
\psi_{ v} : G \times U \to F \qquad \psi_{v}(g,u)=g( v+u)\,.
\end{equation}
The morphism $\psi_{ v}$ induces a $G$-equivariant morphism 
\begin{equation}
\label{regular to ss}
\phi_{ v} :  F \to F^{\mathrm{ss}}  \qquad \phi_{v} (\psi_{v}(g,u))=g( v)
\end{equation}
and its fiber can be identified with the prehomogenous vector space $(Z_G(\Cartan)^0,U)$.

The dual representation $V^*$ is also polar and we have
\begin{equation}
V^*\ = \ U^* \oplus \Cartan^* \oplus \fg\cdot \Cartan^*
\end{equation}
as $Z_G(\Cartan)^0$-representations and as before we have a description for the strongly regular elements analogous to~\eqref{reg decomp}. 

 As in~\cite{G1,G2,GVX1,GVX2} we use (stratified) Morse theory to describe  $\widehat {P_\cF}|_{(V^*)^\reg}$. First, we regard the elements $\xi\in (V^*)^{\reg}$ as linear functions on $V$. Then, arguing as in~\cite{GVX1,GVX2}, the stalk of $\widehat {P_\cF}|_{(V^*)^\reg}$ at the point $\xi$ is given by 
\beqn
(\widehat {P_\cF}|_{(V^*)^\reg})_\xi  \cong H^0 (F, \{ x \in F\; | \; \on{Re} \xi(x) \leq - \xi_0 \}; \cF),
\eeqn
for $\xi_0\gg 0$. The right hand side can be expressed using Morse theory in terms of Picard Lefschetz classes as follows. 

To find the critical points we intersect the characteristic variety $\on{SS}(\cF)$ with the graph of the differential of the  linear form $\xi$, i.e., we determine
\beqn
\{(x,\xi) \mid x\in  V\} \cap \on{SS}(\cF)\,.
\eeqn
By equivariance it suffices to consider elements $\xi\in (U^*)^{\mathrm{r}}\times (\Cartan^*)^\mathrm{rs}$.  Consider the component $T^*_{F^{\mathrm{ss}}} V$ of $\on{SS}(\cF)$. At a point $x\in F^{\mathrm{ss}}$ we have 
\begin{equation*}
(T^*_{F^{\mathrm{ss}}}V)_x = \{\eta\in V^*\mid \eta|_{\Lg. x} =0\} \,.
\end{equation*}
It is not difficult to see that this conormal bundle can intersect the graph of $\xi$ only if $x\in\Cartan^{\mathrm{rs}}$. In that case $\Lg. x= \Lg. \fa$ and we have 
\beqn
(T^*_{F^{\mathrm{ss}}}V)_x\  = \  \fa\oplus U\qquad \text{hence}\qquad T^*_{F^{\mathrm{ss}}}V|_{\Cartan^{\mathrm{rs}}}=(\Cartan\cap F^{\mathrm{ss}})\times( \fa\oplus U)\,.
\eeqn
Now $\Cartan\cap F^{\mathrm{ss}} = W\cdot a_0^{\mathrm{rs}}$ for some element $a_0^{\mathrm{rs}}\in\Cartan^{\mathrm{rs}}$.
Thus we have 
\beqn
\{(x,\xi) \mid x\in V\} \cap T^*_{F^{\mathrm{ss}}} V = \{(x,\xi) \mid x\in V\}\cap (W\cdot a_0^{\mathrm{rs}}\times( \Cartan \oplus U)) = W\cdot a_0^{\mathrm{rs}}\times\{\xi\}.
\eeqn
It is not difficult to verify that this intersection is transverse. Furthermore, the critical points $W\cdot a_0^{\mathrm{rs}}$ are independent of the Morse function  $\xi\in (U^*)^{\rr}\times (\Cartan^*)^{\mathrm{rs}}$. 

We claim that the graph of $\xi$ does not intersect any other conormal bundles $T^*_{F'}V$ if $F'$ is not the semisimple orbit in the fiber $F$. This is easy to see from the decomposition~\eqref{decomp}. One observes that the form $\xi$ is non-trivial on the tangent space to $F'\cap \phi_v^{-1}(v)$.

We can now give a geometric description of $\widehat {P_\cF}|_{(V^*)^\reg}$. Consider the microlocalisation $\mu_{F^{\mathrm{ss}}}(\cF)$. The $\mu_{F^{\mathrm{ss}}}(\cF)$, often also referred to as normal Morse data, is a local system on 
\beqn
\mathring{ T}^*_{F^{\mathrm{ss}}}V \ := \  T^*_{F^{\mathrm{ss}}}V \cap (V\times (V^*)^{\reg})\,.
\eeqn
Using the Picard-Lefschetz cuts, as explained in~\cite{GVX1,GVX2}, we can express the stalk $(\widehat {P_\cF}|_{(V^*)^\reg})_\xi$   in the following terms:
\beq
\label{Morse data}
(\widehat {P_\cF}|_{(V^*)^\reg})_\xi  \ = \ \bigoplus_{w\in W} TM_w \otimes \mu_{F^{\rrs}}(\cF)_\xi\,.
\eeq
Here the $TM_w$ denotes the tangential Morse group at $w\cdot a_0^{\mathrm{rs}}$.  Note that when $\xi$ moves in $(U^*)^{\rr}\times (\Cartan^*)^{\mathrm{rs}}$, the first term in the tensor product in~\eqref{Morse data} does not change, but the second term does depend on $\xi$ on its  $(U^*)^{\rr}$ component. The isomorphism in~\eqref{Morse data} depends on the  Picard-Lefschetz cuts as they vary continuously when we move $\xi$. In  \autoref{nearby} we carry out this analysis in detail in the context of graded Lie algebras following the general strategy of~\cite{GVX2}. 

Finally,  if $\cF$ is irreducible, we have the following generalization of the result of Grinberg in~\cite{G2} :
\beq
\label{IC}
\widehat {P_\cF}  = \  \IC( \widehat {P_\cF}|_{(V^*)^\reg}).
\eeq
In~\cite{GVX3} we prove this result, even more generally for arbitrary $q\in Q$, when the polar representation arises from the graded Lie algebras as explained in  \autoref{graded}. Thus, what needs to be understood is the local system $ \widehat {P_\cF}|_{(V^*)^\reg}$. To explicitly determine $ \widehat {P_\cF}|_{(V^*)^\reg}$ we will later assume, moreover, that $\cF$ is a clean extension from $F^\reg$ .  In this case and under the exactness condition~\eqref{exact hypothesis} the statement~\eqref{IC} follows directly from the adjunction formula~\cite[Theorem 3.3]{LTVX}. This argument is given in~\autoref{Proof of IC}. By making use of~\eqref{regular to ss} we can identify $G$-equivariant local systems on $F^{\reg}$ with $Z_G(\Cartan)$-equivariant local systems on $U$. Let us note again that $(Z_G(\Cartan)^0,U)$ is a prehomogenous vector space and also a rank zero polar representation. We will discuss this situation in the next section. 

In our applications the restriction of $\mu_{F^{\mathrm{ss}}}(\cF)$ to $\mathring{ T}^*_{F^{\mathrm{ss}}}$ will be a local system of rank one. Hence the local system $ \widehat {P_\cF}|_{(V^*)^\reg}$ is of rank $|W|$ just like what we obtain from the invariant system construction. However, these two constructions do not in general coincide. Later in this paper we will show that in some very special cases $ \widehat {P_\cF}|_{(V^*)^\reg}$ coincides with the restriction to ${(V^*)^\reg}$ of an invariant system (see Proposition~\ref{First estimate}).

\subsection{Nearby cycle functor via the ${\mathbf{V}}$-filtration}We conclude this section by recalling the nearby cycle construction for $\cD$-modules in the form utilised in the rest of the paper, i.e., via the ${\mathbf{V}}$-filtration.  
We consider $f:V\to \bC$ where $f$ can be any function on a smooth variety $V$.  Let $\cN$ denote any regular holonomic $\cD_V$-module on $V$.  

The construction of nearby cycles using ${\mathbf{V}}$-filtration requires us to use a smooth divisor, so we first pass to that situation. Consider
\beqn
 V \xrightarrow{i} V\times \bC \xrightarrow{p} \bC \qquad i(v)= (v,f(v)) \ \ \ p(v,a)=a\,
\eeqn
so that $f=p\circ i$. We  form $i_*\cN$ and then take nearby cycles with respect to $p$.  Recall that for the $\mathbf{V}$-filtration the elements in $\cD_V$ are considered to be of degree 0, $t$ to be of degree 1 and $\del_t$ to be of degree $-1$, where $t$ is the coordinate of $\bC$. 
The nearby cycles are given by
\beqn
\psi_f \cN\ = \ \oplus_{-1<\alpha\leq 0} \Gr_{\mathbf{V}}^\alpha i_*\cN\qquad \text{where}\qquad \Gr_{\mathbf{V}}^\alpha i_*\cN = V^{\geq\alpha}i_*\cN/V^{>\alpha}i_*\cN.
\eeqn
The $ \Gr_{\mathbf{V}}^\alpha i_*\cN$ are generalized $\alpha$-eigenspaces for $t\del_t$.
In general the $\alpha$ can be complex numbers but in the cases considered in this paper they are rational numbers. We also recall that the vanishing cycles are given by
\beqn
\phi_f \cN\ = \ \oplus_{-1\leq\alpha< 0} \Gr_{\mathbf{V}}^\alpha i_*\cN\,.
\eeqn

The definitions above are in terms of $i_*\cN$ on $V\times \bC$. We will now express the nearby cycles and vanishing cycles purely in terms of $V$. First observe that we have 
\beqn
i_*\cN \ =\frac{\cO_{V\times \bC}[(t - f)^{-1}]}{\cO_{V \times \bC}} \otimes_{\cO_V} \cN\,.
\eeqn
 For sections $u$ of $\cN$, consider the map 
\begin{equation}
\cD_V[s] f^{s}u \ \ \longrightarrow \ \ \cD_V[\del_t t] \frac {u}{t-f}\,
\end{equation}
which takes
\begin{equation}\label{eqn-hom}
\begin{aligned}
s \ \ &\longmapsto \ \ -\del_t t
\\
f^su  \ \ &\longmapsto \ \ \frac {u}{t-f}\,.
\end{aligned}
\end{equation}
Now consider the Bernstein-Sato b-function of $u$:
\beqn
D f^s u \ = \ b(s)f^{s-1} u  \,.
\eeqn
Note that we use a slightly different normalisation of the b-function as it is more convenient for our purposes and results in simpler formulas.
We then write
\beqn
D f^{s+1} u \ = \ b(s+1)f^{s} u\,.
\eeqn
Noting that 
\beqn
\frac f {t-f} u = \frac t {t-f} u
\eeqn and applying~\eqref{eqn-hom}, the above equation becomes
\begin{equation}
\label{b function for V-filtration}
Dt\,\frac {u}{t-f} \ =\ b(-\del_t t+1)\,\frac {u}{t-f}\,.
\end{equation}
Setting $b_\bfV(t\del_t)=b(-t\del_t) = b(-\del_t t+1)$ we obtain
\beqn
Dt\,\frac {u}{t-f} \ =\ b_\bfV(t\del_t)\,\frac {u}{t-f}\,
\eeqn
i.e., the ${\mathbf{V}}$-filtration b-function equation. Since the $b(s)$ is the b-function of $u$, $b_\bfV(t\del_t)$ is the ${\mathbf{V}}$-filtration b-function for $\frac {u}{t-f}$. Thus, we have $b_\bfV(s)=b(-s)$. 

Assume now that the section $u$ generates $\cN$  and let $b_u$ denote its ${\mathbf{V}}$-filtration b-function. We write $b_0$ for the polynomial obtained from $b_u$ by shifting all of its zeros to the interval $(-1,0]$. It is a polynomial of the same degree as $b_u$. By considerations in~\cite[Section 4]{maisonbe-mebkhout}, for example, we see that 
\begin{proposition}
\label{monodromy from b-function}
The $b_0(s)$ is the minimal polynomial for the operator $t\del_t$ on $\psi_f \cN$. 
\end{proposition}

\section{Invariant systems in rank zero}
\label{zero}

In this section we consider invariant systems when the polar representation has rank zero. We will make use of results of~\cite{Sa,SK}  and also use~\cite{K} as a convenient reference.  The rank zero case amounts to the case of  {\em finite} prehomogenous vector spaces $(G,V)$, i.e., the vector space $V$ consists of finitely many $G$-orbits. We continue in the general notation of the previous section. We will assume in addition that $(G,V)$ is {\em regular}, that is, the open $G$-orbit $V^\reg$ in $V$ is affine.  The situation we have in mind is the space $\Lu_1$ with the $Z_K(\fa)^0$-action of~\autoref{graded}. In our application in~\autoref{cuspidal sheaves} of the theory developed here, the stabilizer of a generic point is finite.

\subsection{Fundamental semi-invariants and the invariant systems $\cM_{s_\bullet}$}We write $\rho:G \to GL(V)$ for the representation.  Recall that a function $h\in \cO(V)$ is called a semi-invariant for $(G,V)$ if there is a character $\chi:G \to \bC^*$ such that $h(gv) = \chi(g)h(v)$ for all $v\in V$, $g\in G$. We write $\on{SI}(G,V)$ for the set of semi-invariants. They form a sublattice in  $X^*(G)$. We write $f_1, \dots f_k$ for the fundamental semi-invariants  and $\chi_i:G \to \bC^*$, $i=1,\ldots,k$, for the corresponding characters. The fundamental semi-invariants are irreducible. 

Let  $G_1=G_{der}G_v$ where $v$ is an arbitrary element in $V^\reg$. The group $G_1$ is independent of the choice of $v$ and the semi-invariants factor through the torus $G/G_1$.

Consider a tuple $s_\bullet=(s_1,\dots, s_{k})$ of elements in $\bQ$ such that 
\beq\label{eqn-chis}
\chi=\chi_{s_\bullet}=\prod_{i=1}^{k} \chi_i^{s_i} \in X^*(G)\,.
\eeq
As was pointed out in the last section, under this condition the invariant system 
\beq
\cM =\cM_ {s_\bullet}:=\cD_V / \langle  \{d\rho(\xi)- \sum_{i=1}^{k} s_id\chi_i(\xi) \mid \xi\in \fg \}\rangle
\eeq
is $G$-equivariant.

As $f_i$ is a semi-invariant we see that $(\xi.f_i)(v) = - d\chi_i(\xi)f_i(v)$. Thus, we conclude that 
\beq
\cM|_{V^\reg} \cong \cD_{V^\reg} f_1^{-s_1}\dots f_k^{-s_k}\,.
\eeq
Note that under the deRham functor the local system corresponding to the $\cD_{V^\reg}$-module $\cD_{V^\reg} f_1^{-s_1}\dots f_k^{-s_k}$ is given by the monodromy of the function $f_1^{s_1}\dots f_k^{s_k}$, which is exactly the monodromy of the local system associated to the character $\chi|_{G_v/(G_v)^0}$; here $v$ is any element in $V^\reg$. We write $$\cL(s_1,\dots,s_{k})=\cL_\chi=\cL_{\chi_{s\bullet}}$$ for the resulting rank one local system on $V^{\reg}$. That is, $\cM|_{V^\reg}$ corresponds to $\cL_\chi$ under the deRham functor.

Making use of~\cite[Proposition 2.11]{K} and observing that a rank one local system on $V^\reg$ is determined by its monodromies around the divisors $f_i=0$, i.e., is given by the monodromy of a function $f_1^{s_1}\dots f_k^{s_k}$ for some $s_i\in\bC$,  we conclude that
\begin{enumerate}
\item[(i)] all $G$-equivariant rank one local systems on $V^\reg$ arise from such characters $\chi$ as in~\eqref{eqn-chis}
\item[(ii)] two such characters give rise to the same local system precisely when the  corresponding $(s_1,\dots,s_{k})$ coincide modulo integers.
\end{enumerate}

\subsection{Fourier transform of $\cM_{s_\bullet}$}Let us recall the Fourier transform in the language of $\cD$-modules. 
 Write $\rho^*:G \to GL(V^*)$ for the action of $G$ on the dual representation $V^*$. The group $G_1^*=G_{der}G_{v^*}=G_1$, with $v^*\in (V^*)^\reg$. The fundamental semi-invariants $f_i^*\in\cO(V^*)$ are associated to characters $\chi_i^*=\chi_i^{-1}\in X^*(G)$. The action of the group $G$ on $V$ gives rise to a map $d\rho:\fg \to V \otimes V^*$. We regard $V \otimes V^*$ as elements in $\cD_V$ in the usual way, i.e., to an element $(D,x)\in V \otimes V^*$ we associate the operator $xD$ and in this manner we get a map $d\rho: \fg \to \cD_V$. Similarly, we have a map $$d\rho^*: \fg \to \cD_{V^*}.$$ For concreteness we fix dual bases $\del_i$ of $V$ and $x_i$ of $V^*$. We will also write $\del_i= x_i^*$ and $x_i=\del_i^*$ for the basis elements. 

 By definition, the Fourier transform $\widehat\cM$ of the $\cD_V$-module $\cM$ is the module $\cM$ regarded as a $\cD_{V^*}$-module $\cM$ via the identification of the two algebras of differential operators as follows
\beqn
\widehat{} :\cD_V \xrightarrow{\ \ \cong \ \ }  \cD_{V^*} \qquad    \hat x_i = - \del_i^* \ \ \  \hat \del_i = x_i^*\,.
\eeqn
A short calculation shows that
\beqn
\widehat {d\rho(\xi)} = d\rho^*(\xi)  + \on{Tr}(d\rho(\xi))\,.
\eeqn
We can also regard the  term $ \on{Tr}(d\rho(\xi))$ as the differential of the character $\det(\rho): G \to \bC^*$. Thus we see that the Fourier transform of 
\begin{subequations}
\beq
\cM = \cD_V / \langle  \{d\rho(\xi)- \sum_{i=1}^{k} s_id\chi_i(\xi) \mid \xi\in \fg \}\rangle
\eeq
is
\beq
\widehat\cM = \cD_{V^*} / \langle  \{d\rho^*(\xi)+ \sum_{i=1}^{k} s_id\chi_i^*(\xi) - \on{Tr}(d\rho^*(\xi)) \mid \xi\in \fg \}\rangle\,.
\eeq
\end{subequations}
 Recall that $\det(\rho)^2$ vanishes on $G_1$ and hence is a semi-invariant. Thus, there exist rational numbers $s^*_i$ such that the module $\hat\cM$ can  be written in the following manner 
\beq
\widehat\cM = \cD_{V^*} / \langle   \{d\rho^*(\xi)- \sum_{i=1}^{k} s^*_i d\chi_i^*(\xi) \mid \xi\in \fg \}\rangle\,.
\eeq
In particular, $\widehat\cM$ is also an invariant system. 

\subsection{Clean invariant systems}The constructions in this section will be used in very specific settings later in the paper, i.e., in the cases when the invariant system $\cM_{s_\bullet}$ is clean. In this paper a key role is played by b-functions. We recall here how to characterise cleanness in terms of the b-function. Let us write, just in this subsection, 
\beqn
u=  f_1^{-s_1}\dots f_k^{-s_k}\text{ and }f=f_1 \dots f_k.
\eeqn 
Consider the $f_i^*$ as (constant coefficient) differential operators and denote them by $D_i$. Let $D=D_1 \dots D_k$. By the semi-invariance property we know that the b-function is given by 
\beq
\label{b-fcn ph}
D f^s u \ = \ b(s)f^{s-1} u  \,.
\eeq
Recall that we use a slightly different normalisation of the b-function.  

The local system $\cL_\chi=\cL(s_1,\ldots,s_k)$ only depends on the $s_i$ modulo $\bZ$. However, that is not the case for the invariant system $\cM=\cM_{s_\bullet}$. Let $j: V^\reg \to V$ be the inclusion. The b-function lemma implies the following well-known proposition:
\begin{proposition}\label{clean criterion}
The invariant system $\cM_ {s_\bullet}$ is is a clean extension if and only if $b(s)$ in~\eqref{b-fcn ph} has no integer zeros. \end{proposition}

\section{Invariant systems in rank one}\label{sec-rank1g}

In this section we consider invariant systems for rank one polar representations.

Let $(G,V)$ be a {\em visible} polar representation of rank one and we write $\rho:G \to GL(V)$ for the representation.    The $G\times\bC^*$-action makes $V$ into a prehomogenous vector space with $V^\reg$ as the open dense orbit. We assume that the pair $(G\times\bC^*,V)$ constitutes a regular finite prehomogenous vector space. This is the case if the pair arises from a rank one graded Lie algebra.

Let us write $\bC[V]^G = \bC[f]$. 
Let  $f_0, \dots, f_k$  be the fundamental semi-invariants which appear as factors of  the invariant $f$. Thus 
\beq\label{eqn-inv}
f= f_0^{n_0} \dots f_k^{n_k}
\eeq for some positive integers $n_i$. Let $f_{k+1}, \dots, f_{\ke}$,  be the remaining fundamental semi-invariants (where $e\geq 0$). We write  $\chi_i:G \to \bC^*$, $i=0,\ldots,\ke$, for the corresponding characters where we have omitted the $\bC^*$-factor which keeps track of the degree of homogeneity of the semi-invariants. After this restriction the characters $\chi_i$ are no longer independent but satisfy the relation 
\beq
\chi_0^{n_0} \dots \chi_k^{n_k} \ = \ 1\,.
\eeq
The semi-invariants $\on{SI}(G,V)$ form a sub-lattice of all characters $X^*(G)$. Let $G_1=G_{der}G_v$ where $v$ is an arbitrary element in $V^\reg$. The group $G_1$ is independent of the choice of $v$ and the semi-invariants factor through the torus $G/G_1$.

\subsection{Kostant slice}We will now make a further assumption that there exists a Kostant slice $\Ls$ for $$f: V \to V\inv G = \bC.$$ This is equivalent to  at least one of the $n_i$ in~\eqref{eqn-inv} being equal to one. We will  from now on assume that $n_0=1$. Thus, the Kostant slice $\Ls$ is a germ of an analytic variety which meets the zero locus of $f_0$ transversally but does not meet the zero locus of any other semi-invariants. We make use of the Kostant slice on the dual side. 
To that end consider the dual representation $(G,V^*)$. We have 
\beqn
\chi^*_0 (\chi_1^*)^{n_1} \dots (\chi_k^*)^{n_k} \ = \ 1\,.
\eeqn
Thus, the invariant is
\beqn
f^* \ = \ f^*_0 (f_1^*)^{n_1} \dots (f_k^*)^{n_k}\,.
\eeqn
We will consider the Kostant slice $\Ls^*$ which meets the zero locus of $f^*_0$ transversally but not the zero locus of any other semi-invariants. Fix $v^*\in (V^*)^\reg$. We have the following exact sequence:
\beqn
1 \to G_{v^*}/ G_{v^*}^0 \to \pi_1^G\left((V^*)^\reg, v^*\right) \to \pi_1((V^*)^\reg\inv G) =\bZ \to 1\,.
\eeqn
The Kostant slice splits this exact sequence by identifying $\Ls^*\cap (V^*)^\reg$ with $(V^*)^\reg\inv G$. Then   $\pi_1^G((V^*)^\reg, v^*)\cong G_{v^*}/ G_{v^*}^0 \rtimes \bZ$. Completely analogously, we have the exact sequence 
\beqn
1 \to G_{v}/ G_{v}^0 \to \pi_1^G(V^\reg, v) \to \pi_1(V^\reg\inv G) =\bZ \to 1
\eeqn
which splits similarly and we have $\pi_1^G(V^\reg, v)\cong G_{v}/ G_{v}^0 \rtimes \bZ$. As a matter of fact the spaces $V^\reg$ and $(V^*)^\reg$ can be identified as explained in~\cite[Theorem 2.16]{K}. This identifies the two exact sequences above.

\subsection{The nearby cycles}We will consider all rank one $G$-equivariant $\cD$-modules on $V^\reg$ such that the holonomy around the divisor given by $f_0=0$ is trivial. Such a $G$-equivariant $\cD$-module can be restricted to a   $G$-equivariant $\cD$-modules on $f^{-1}(a)^\reg = f^{-1}(a) \cap V^\reg$ where $a$ is any non-zero element in $\bC$. We can understand these $\cD$-modules in the same manner as in the previous section. 
Consider a tuple $s_\bullet=(s_1,\dots, s_{\ke})$ of elements in $\bQ$ such that 
\beq\label{eqn-tuples}
\chi=\chi_{s_\bullet}=\prod_{i=1}^{\ke} \chi_i^{s_i} \in X^*(G)\,.
\eeq  
The character $\chi$ gives rise to  a character of $G_v/(G_v)^0$ and hence a local system on $f^{-1}(a)^\reg$ which is the restriction of a local system on $V^\reg$. Its monodromy on $V^\reg$ is given by the holonomy of the multivalued function $f_1^{s_1}\dots f_{\ke}^{s_{\ke}}$. Via the deRham functor this local system coincides with the $\cD$-module $\cD_{V^\reg} f_1^{-s_1}\dots f_{\ke}^{-s_{\ke}}$. 

We note that the above construction does not necessarily yield all rank one $G$-equivariant local system on $f^{-1}(a)^\reg$. However, we obtain all such local systems if the action of $\pi_1(V^\reg\inv G) =\bZ$ is trivial on $G_{v}/ G_{v}^0$ in which case we have  $\pi_1^G(V^\reg, v)\cong G_{v}/ G_{v}^0 \times \bZ$.

Set 
\beqn
\cN \ =\ \cD_{V}u,\ u=  f_1^{-s_1}\dots f_{\ke}^{-s_{\ke}} \,.
\eeqn
Let 
\beqn
P_\chi  = \psi_f \cN\,.
\eeqn
It is a $\cD$-module supported on the null-cone $f^{-1}(0)$.

We will now impose the following condition on $\cN$:
\beq
\text{For $a\in \bC^* = \Cartan^{rs}/W$ the $\cD$-module $\cN|_{f^{-1}(a)}$ is a clean extension from $f^{-1}(a)^\reg$.}
\eeq
We can determine the clean extensions as follows. We analyse the fibre $F= f^{-1}(a)$ as in~\eqref{regular to ss}. In particular, recall that the fiber of $\phi_v$ can be identified with the prehomogenous vector space $(Z_G(\fa)^0, U)$. 

As the $f_i$ are $G\times \bC^*$ semi-invariants, they give rise to non-trivial semi-invariants $\tilde f_i$  of $Z_G(\fa)^0$ on $U$. We now work on $U$ and write $\tilde u= \tilde f_{k+1}^{-s_{k+1}}\dots \tilde f_{\ke}^{-s_{\ke}}$, $\tilde f= \tilde f_{k+1}\dots \tilde f_{\ke}$ and $\tilde D$ for the corresponding differential operator. Again by semi-invariance we have 
\beq\label{}
\tilde D \tilde f^s \tilde u \ = \ \tilde b(s)\tilde f^{s-1} \tilde u
\eeq for the b-function. 
Then, making use of Proposition~\ref{clean criterion} we have
\beq
\label{clean criterion rk1}
\text{$\cN|_{f^{-1}(a)}$ is a clean extension from $f^{-1}(a)^\reg$ if and only if $\tilde b(s)$ has no integer zeros.}
\eeq

\subsection{The Fourier transform}Our goal is to determine the local system $\widehat P_\chi|_{(V^*)^\reg}$ as a representation of $G_{v^*}/ G_{v^*}^0 \rtimes \bZ$.
In what follows we need some control on the endomorphisms of $P_\chi $. We can achieve that by assuming that~\eqref{IC} holds, that is,  
\beq
\label{IC special}
\widehat P_\chi  = \  \IC( \widehat P_\chi|_{(V^*)^\reg})\,.
\eeq
Recall that, as explained in~\autoref{invariant},~\eqref{IC special} holds by~\cite{GVX3} when the polar representation arises from a graded Lie algebra.  As the proof is not yet available in the literature, we will give an alternative argument in \autoref{Proof of IC} which suffices for the applications to cuspidal character sheaves in \autoref{cuspidal sheaves}.  
\begin{remark}
 We also suspect that~\eqref{IC} holds under the following additional hypothesis on the polar representation:
\beq
\label{strongly ruled}
\text{There exists a complement $\fn$  to $Z_\fg(\fa)$ in $\fg$  such that $\fn\cdot U \subset \fg\cdot \fa$.}
\eeq
\end{remark}

We will first argue that $\widehat P_\chi|_{(V^*)^\reg}$ is a local system of rank $|W|$. Appealing to~\eqref{Morse data} we need to show that the microlocalization $\mu_{F^{\mathrm{ss}}}(\cN|_{f^{-1}(a)})$ is of rank one. Using notation from \autoref{invariant} we restrict $\cN|_{f^{-1}(a)}$ further to $U$. Then $\cN|_U$ is a $Z_G(\fa)^0$-equivariant clean sheaf on $U$. Recall that $U$ is a prehomogenous vector space for $Z_G(\fa)^0$. We can now apply the considerations from \autoref{zero} and conclude that the Fourier transform of $\cN|_U$ is of rank one. Arguing as in~\cite{GVX2} we can further conclude that
\beq
\label{cyclic}
\text{the local system $\widehat P_\chi|_{(V^*)^\reg}$ as a representation of $\bZ$ is cyclic,}
\eeq
 i.e., it is given by a regular matrix.

\subsection{Comparing the nearby cycle sheaf and the invariant system} \label{ssec-comparison}We now  analyse the relationship between $P_\chi$ and the invariant system
\beq
\label{rk1 M}
\cM_\chi \ =\  \cD_V / \langle f, \{d\rho(\xi)- \sum_{i=1}^{l+k} s_id\chi_i(\xi) \mid \xi\in \fg \}\rangle\,.
\eeq 
The system $\cM_\chi$ is $G$-equivariant and is supported on the null-cone by construction. 

The Fourier transform of 
\begin{subequations}
\label{invt and ft}
\beq
\cM_\chi = \cD_V / \langle f, \{d\rho(\xi)- \sum_{i=1}^{l+k} s_id\chi_i(\xi) \mid \xi\in \fg \}\rangle
\eeq
is
\beq
\widehat\cM_\chi = \cD_{V^*} / \langle \hat f = D^*, \{d\rho^*(\xi)+ \sum_{i=1}^{l+k} s_id\chi_i^*(\xi) - \on{Tr}(d\rho^*(\xi)) \mid \xi\in \fg \}\rangle\,.
\eeq
\end{subequations}
In the above $f$ is the generating invariant in $(\on{Sym}V^*)^G$ and $\hat f = D^*$ is the same generating invariant but viewed as a differential operator on $V^*$. Recall that  $\det(\rho)^2$  is a semi-invariant. Hence we can write 
\beq
\label{expression for det}
(\on{det}\rho)^2=\prod_{i=0}^{\ke}\chi_i^{2\kappa_i}
\eeq
where the $\kappa_i$ are half integers. This expression is unique if we insist that it holds for the $G\times \bC^*$-action. In what follows we assume that we have made this choice. It will play an important role when we discuss the functional equation for b-functions.

Recall also that $\chi_i^*=\chi_i^{-1}$. Then $ \on{Tr}(d\rho^*(\xi))= \sum_{i=0}^{\ke}\kappa_i d\chi_i^*(\xi)$. As $\chi_0 = \chi_1^{-n_1} \dots \chi_k^{-n_k}$ we see that
\beq
 \on{Tr}(d\rho^*(\xi))= \sum_{i=1}^{k}(\kappa_i- n_i\kappa_0)d\chi_i^*(\xi) +  \sum_{i=k+1}^{\ke}\kappa_i d\chi_i^*(\xi)\,.
\eeq
Thus we can write $\widehat\cM_\chi$ in the following manner 
\beq
\widehat\cM_\chi = \cD_{V^*} / \langle  D^*, \{d\rho^*(\xi)- \sum_{i=1}^{l+k} s^*_i d\chi_i^*(\xi) \mid \xi\in \fg \}\rangle\,
\eeq
where
\beq
\label{s* in terms of s}
s^*_i= 
\begin{cases} -s_i-n_i\kappa_0+\kappa_i, &\text{if} \ \  1\leq i \leq k  
\\
-s_i +\kappa_i &\text{if} \ \  k+1\leq i \leq k +l\,.
\end{cases}
\eeq

Let $D=f^*$ be the generating invariant of $\on{Sym}(V)^G$ viewed as a differential operator on $V$. Consider the b-function
\beq
\label{b-fcn}
D f^s u \ = \  b(s)f^{s-1} u  \,,\quad u=  f_1^{-s_1}\dots f_{\ke}^{-s_{\ke}}.
\eeq
As before we use a slightly different normalisation of the b-function.  
By $G$-equivariance we know that the b-function equation takes this form. 

Let us recall the discussion about the nearby cycle functor via the ${\mathbf{V}}$-filtration. In particular, recall that the ${\mathbf{V}}$-filtration b-function $b_\bfV$ of $\frac {u}{t-f}$ is given by $b_\bfV(s)= b(-s)$ (up to a non-zero constant). We will follow the notation of \autoref{invariant}.

\begin{proposition}
\label{First estimate}
Assume that the zeros of the b-function $b(s)$ in~\eqref{b-fcn} are located in within a half-open interval $J$ of length 1 containing zero.  Then $P_\chi $ is a quotient of $\cM_\chi$. In particular, $\widehat P_\chi|_{(V^*)^\reg} = \widehat\cM_\chi|_{(V^*)^\reg}$.
\end{proposition}

\begin{proof}

Let us write (up to constant) $b(s) = \prod_{i=1}^p (s+\alpha_i)^{d_i}$ so that $b_\bfV(s) = \prod_{i=1}^p (s-\alpha_i)^{d_i}$. Note that $\deg b(s) = \sum d_i = n$ where $n= \deg f = \deg D$. We have
\beq
P_\chi  = \psi_f \cN = \oplus_{i=1}^p \Gr_{\mathbf{V}}^{\alpha_i}i_* \cN\,.
\eeq
The $ \Gr_{\mathbf{V}}^{\alpha_i}i_*\cN$ are generalized $\alpha_i$-eigenspaces for $t\del_t$. By our assumption of the roots of the b-function, the polynomial $b_\bfV(s)$ is the minimal polynomial for $t\del_t$ acting on $P_\chi$.
Assume that we have numbered the roots in increasing order so that $\alpha_1$ is the smallest. We then have the following elements in  $V^{\alpha_i}i_*\cN$ which give us non-zero elements $\bar u_i\in \Gr_{\mathbf{V}}^{\alpha_i}i_*\cN$:
 \beq
 u_i= \prod_{e=1}^{i-1} (t\del_t-\alpha_e)^{d_e}  \frac u {t-f}  \in V^{\alpha_i}\cN\,.
 \eeq
 Now,
 \beqn
fu_i=  \prod_{e=1}^{i-1} (t\del_t-\alpha_e)^{d_e}\frac f {t-f} u =  \prod_{e=1}^{i-1} (t\del_t-\alpha_e)^{d_e}\frac t {t-f} u \in V^{\alpha_i+1}i_*\cN\,.
\eeqn
By making use of our hypothesis about the zeros of $b$ we conclude that $f\bar u_i = 0$. As $f$ is an invariant, we have $d\rho(\xi)f=0$ for $\xi\in\fg$. Recalling that  $(\xi.f_i)(v) = - d\chi_i(\xi)f_i(v)$ and by definition of $\cN$ we see that  $(d\rho(\xi)- \sum_{i=1}^{l+k} s_id\chi_i(\xi))u_i=0$ for $\xi\in\fg$. Thus we get a morphism
\beq
\cM_\chi \to P_\chi= \oplus_{i=1}^p \Gr_{\mathbf{V}}^{\alpha_i}i_* \cN \qquad 1\mapsto \bar u = (\bar u_1, \dots, \bar u_p)
\label{M to P}
\eeq
and so a morphism 
\beq
\widehat \cM_\chi \to \widehat P_\chi  = \  \IC( \widehat P_\chi|_{(V^*)^\reg})\,.
\eeq
To show that~\eqref{M to P} is a surjection we  argue that $\bar u$ generates $ P_\chi$ as a $\cD$-module. To do so we first argue that the endomorphism $t\del_t$ has minimal polynomial $(s-\alpha_i)^{d_i}$ on the section $\bar u_i\in \Gr_{\mathbf{V}}^{\alpha_i}i_*\cN$. But, the polynomial $b_\bfV(t\del_t)$ is the minimal polynomial which moves the zeros $\alpha_i$ out of the half-open interval $J$. This implies that the minimal polynomial moving the zeros of $u_i$ out of the half-open interval $J$ is $\prod_{e=i}^{p} (t\del_t-\alpha_e)^{d_e}$. Thus, the minimal polynomial for each $\bar u_i$ is $(t\del_t-\alpha_i)^{d_i}$. Hence the minimal polynomial for $\bar u$ is $b_\bfV(t\del_t)$. Now consider the Fourier transform $\widehat {\cD_V \bar u}$. We have 
\beqn
\widehat {\cD_V \bar u}   \subset \widehat P_\chi = \IC( \widehat P_\chi|_{(V^*)^\reg})\,.
\eeqn
Therefore we conclude that 
\beq
\widehat {\cD_V \bar u}= \IC( \widehat {\cD_V \bar u}|_{(V^*)^\reg})\,.
\eeq
The operator $t\del_t$ induces an endomorphism of $\widehat {\cD_V \bar u}$ and because $\widehat {\cD_V \bar u}$ is an $\IC$-extension this endomorphism amounts to an endomorphism of the local system  $\widehat {\cD_V \bar u}|_{(V^*)^\reg} $. We had just argued that the minimal polynomial of this endomorphism is of degree $n=|W|$ and hence the local system has to have this rank. This is also the rank of  $\widehat P_\chi|_{(V^*)^\reg}$ and thus we have 
$\widehat {\cD_V \bar u} =  \widehat P_\chi$ which implies that $\cD_V \bar u= P_\chi$.
\end{proof}

\subsection{The local system $ \widehat\cM_\chi|_{(V^*)^\reg}$}We will now determine the local system $ \widehat\cM_\chi|_{(V^*)^\reg}$.
Recall that $(\xi.f^*_i)(v) = - d\chi^*_i(\xi)f^*_i(v)$ and that 
\beqn
\widehat\cM_\chi = \cD_{V^*} / \langle  D^*, \{d\rho^*(\xi)- \sum_{i=1}^{l+k} s^*_i d\chi_i^*(\xi) \mid \xi\in \fg \}\rangle\,.
\eeqn
Let us look at solutions to this system on $(V^*)^\reg$. By $G$-equivariance, the solutions on  $(V^*)^\reg$  are given as follows:
\beqn
\Hom(\widehat\cM_\chi,\cO_{V^*}^{an}) = \{F(f^*) (f_1^*)^{-s^*_1}\dots (f^*_{\ke})^{-s^*_{\ke}}\mid DF(f^*) (f_1^*)^{-s^*_1}\dots (f^*_{\ke})^{-s^*_{\ke}}=0 \}\,.
\eeqn
Let $u^*= (f_1^*)^{-s^*_1}\dots (f^*_{\ke})^{-s^*_{\ke}}$ and  consider the b-function
\beqn
D^* (f^*)^s u^* \ = \  b^*(s)(f^*)^{s-1} u^*. 
\eeqn
Then we have $  f^*D^* (f^*)^s u^* = b^*(s)(f^*)^{s} u^*$.
From this we conclude that the equation $D^*F(f^*)u^*=0$ is equivalent to the following equation for $F(t)$:
\beq
\label{de}
b^*(t\partial_t) F(t) \ = \ 0\,.
\eeq
Let us factor the b-function (up to a nonzero constant)
\beqn
b^*(s) = \prod_{i=1}^p (s-\alpha_i^*)^{d_i}\,.
\eeqn
Then the solution to~\eqref{de} is given explicitly by the formula
\beqn
\sum_{i=1}^p \sum_{h=0}^{d_i-1}  C_{i,h}(\log t)^h t^{\alpha^*_i}.
\eeqn
Assume that the zeros $\alpha_i^*$ lie in a half-open interval of length 1. Then the minimal polynomial of the  monodromy of the local system $\Hom(\widehat\cM_\chi,\cO_{V^*}^{an})|_{(V^*)^\reg}$  coincides with its characteristic polynomial and is given by $\prod (x-\exp(2\pi i \alpha^*_i))^{d_i}$. The monodromy of $ \widehat\cM_\chi|_{(V^*)^\reg}$ via the deRham functor is then given by the polynomial $\prod (x-\exp(-2\pi i \alpha^*_i))^{d_i}$. The representation of $G_{v^*}/ G_{v^*}^0$ is given by the monodromy of the function $(f_1^*)^{s^*_1}\dots (f^*_{\ke})^{s^*_{\ke}}$, i.e., by the character 
\beq
\bega
(\chi_1^*)^{s^*_1}\dots (\chi^*_{\ke})^{s^*_{\ke}} = \prod_{i=1}^k \chi_i^{s_i+n_i\kappa_0-\kappa_i} \prod_{i=k+1}^{\ke} \chi_i^{s_i -\kappa_i } = \chi (\on{det}\rho)^{-1}=\chi\on{det}\rho\,.
\eega
\eeq

We will now explain a functional equation for b-functions which allows us to compare the zeros of the b-functions $b(s)$ and $b^*(s)$.

We write $\underline{\kappa}=(\kappa_0,\ldots,\kappa_{\ke})$ and 
\beqn
\text{$\underline{m}=(1,n_1,\ldots,n_k,0,\ldots,0)$ and $\underline{t}=(0,s_1,\ldots,s_{\ke})$.}
\eeqn 
In the discussion below it is important that we regard $V$ as a $G\times \bC^*$ prehomogenous vector space. In particular it is important that we have chosen the $\underline{\kappa}$ taking into account also the $\bC^*$-action in equation~\eqref{expression for det}.

As before we write $u=\prod_{j=1}^{\ke}f_i^{-s_i}$. 
We have (see~\cite[Theorem 1.3.5]{U2},~\cite[Theorem 2]{Sa})
\beqn
D f^s u=b(s)f^{s-1}u
\eeqn
where, up to a non-zero constant, 
\beqn
b(s)=\prod_{j=1}^N\prod_{\nu=0}^{\gamma_j(\underline{m})-1}\prod_{r=1}^{\mu_j}(s-\frac{\gamma_j(\underline{m})+\gamma_j(\underline{t})-\alpha_{j,r}-\nu}{\gamma_j(\underline{m})})
\eeqn
for some integral linear forms  $\gamma_j(x_0,x_1,\ldots,x_{\ke})=\sum_{i=0}^{\ke}\gamma_{ji}x_i$, for some $N$, $\mu_j$, and rational numbers $\alpha_{j,r}$.

 We write $\underline{t}^*=(0,s_1^*,\ldots,s_{\ke}^*)$ and then by~\eqref{s* in terms of s} $\underline{t}^*=-\underline{t}+\underline{\kappa}-\kappa_0\underline{m}$. 
Recalling that $u^*=\prod_{j=1}^{\ke}(f_i^*)^{-s^*_i}$ we have 
\beqn
D^*(f^*)^s u^*=b^*(s)(f^*)^{s-1} u^*
\eeqn
where, up to a non-zero constant, we have 
\bern
b^*(s)=\prod_{j=1}^N\prod_{\nu=0}^{\gamma_j(\underline{m})-1}\prod_{r=1}^{\mu_j}\left(s-\frac{\gamma_j(\underline{m})+\gamma_j(\underline{t}^*)-\alpha_{j,r}-\nu}{\gamma_j(\underline{m})}\right)\\
=\prod_{j=1}^N\prod_{\nu=0}^{\gamma_j(\underline{m})-1}\prod_{r=1}^{\mu_j}(s-(1-\kappa_0)-\frac{\gamma_j(\underline{\kappa})-\gamma_j(\underline{t})-\alpha_{j,r}-\nu}{\gamma_j(\underline{m})}).
\eern
By~\cite[Corollary 1.3.9]{U2}, for each $j=1,\ldots,N$, we have
\beqn
\{\alpha_{j,r},\,1\leq r\leq\mu_j\}=\{\gamma_j(\underline{\kappa})+1-\alpha_{j,r},\,1\leq r\leq \mu_j\}.
\eeqn
It follows that we can arrange the roots $\alpha_i$ of $b(s)$ and the roots $\alpha_i^*$ of $b^*(s)$ such that
\beqn
\alpha_i+\alpha_i^*=1-\kappa_0.
\eeqn
As $\kappa_0\in\frac{1}{2}\bZ$, this implies that
\beq
\label{zeros relation}
\on{exp}(2\pi\mathbf{i}\alpha_i)\on{exp}(2\pi\mathbf{i}\alpha_i^*)=\on{exp}(2\pi\mathbf{i}\kappa_0)\in\{\pm 1\}\,.
\eeq

As $\pi_1^G((V^*)^\reg, v^*)\cong G_{v^*}/ G_{v^*}^0 \rtimes \bZ$ we can regard ${\widehat \cM_\chi}|_{(V^*)^\reg}$ as a representation of $G_{v^*}/ G_{v^*}^0 \rtimes \bZ$. Thus we obtain a canonical element in $\pi_1^G((V^*)^\reg, v^*)$ corresponding to the element $1\in \bZ$ and we write $x$ for its action on the fiber $(\widehat \cM_\chi)_{v^*}$ and call it the monodromy. 
We summarize the discussion above as a theorem:

\begin{theorem}\label{thm-ft}
Assume that the zeros of $b(s)$ in~\eqref{b-fcn} lie in a half-open interval of length one containing zero.  We write  $b(s)=\prod (s-\alpha_i)^{d_i}$. Then $\widehat P_\chi|_{(V^*)^\reg} \cong\widehat \cM_\chi|_{(V^*)^\reg}$ and, as a representation of $G^{v*}/(G^{v*})^0 \rtimes \bZ$, is given by 
\beqn
\widehat \cM_\chi|_{(V^*)^{\reg}} \ = \ \cL_{\chi\on{det}(\rho)} \otimes \bC[x]/ R_\chi(x), \qquad R_\chi(x) = \prod (x-\exp(2\pi i(\alpha_i+\kappa_0)))^{d_i}\,.
\eeqn
\end{theorem}
\begin{proof}
By~\eqref{zeros relation} we have $\on{exp}(2\pi\mathbf{i}\alpha_i^*)=\on{exp}(-2\pi\mathbf{i}\alpha_i)\on{exp}(2\pi\mathbf{i}\kappa_0)$ and hence
\beq
\prod (x-\exp(-2\pi i \alpha^*_i))^{d_i} = \prod (x- \on{exp}(2\pi\mathbf{i}\alpha_i)\on{exp}(2\pi\mathbf{i}\kappa_0))^{d_i}.
\eeq
\end{proof}
\begin{remark}
Note that in the theorem above the character $\chi\on{det}$ is invariant under the action of $\bZ$. 
\end{remark}

\subsection{Example}

We demonstrate with a simple example the importance of the assumption on the zeros of the b-function. In the example below there are roots of the b-function that do differ by an integer. 

\label{example}
Consider the case of $G=SO_n$ acting on $V=\bC^n$ preserving the quadratic form $f(x) = x_1^2 +\dots x_n^2$. We assume that $n$ is even and $n>2$. In this case all semi-invariants are invariants. Thus the invariant systems are 
\beq
\cM = \cD_V / \langle f, \{d\rho(\xi) \mid \xi\in \fg \}\rangle \ \text{ and }\ \widehat\cM = \cD_{V^*} / \langle \{\Delta, d\rho^*(\xi) \mid \xi\in \fg \}\rangle\,,
\eeq 
where $\Delta$ denotes the Laplacian. 
As an easy calculation shows we have
\beqn
\Delta f^s = 4s(s-1+\frac n 2)f^{s-1},
\eeqn
so that the monic b-function is $b(s)=s(s-1+\frac n 2)$. Thus we have $F(t)=C_1 + C_2t^{-\frac n 2+1}$ which has no monodromy. However, the nearby cycle construction gives us a local system with a non-trivial Jordan form, see, for example~\cite[Example 3.7]{G2}. To get the ``correct" $\cD$-module we have to adjust $\Delta$ with a multiple of the Euler vector field $\theta$. It is  easy to see that the solutions of the $\cD$-module 
\beq
\cD_{V^*}/\langle \{ f\Delta+(2-n)\theta,d\rho^*(\xi) \mid \xi\in\fg \}\rangle
\eeq
are given by $C_1  + C_2 \log f$ and thus have the correct monodromy.

\begin{remark}
We do not have a general guess for an explicit $\cD$-module corresponding to the nearby cycle sheaf. 
\end{remark}

\section{Graded Lie algebras and their equivariant geometry}
\label{graded}

In this section we recall basic facts about graded Lie algebras and explain how the theory in~\cite{GVX2} can be adapted to the setup of this paper.

\subsection{Graded Lie algebras and character sheaves}We follow the definitions and conventions from~\cite{VX2}. Let $G$ be a connected reductive algebraic group over the complex numbers and let $\theta: G \to G$ be an automorphism of $G$ of finite order $m$. We write: 
\beqn
(G^\theta)^0 = K.
\eeqn
Let $\Lg = \on{Lie} G$. Then we have: 
\beqn
\Lg = \bigoplus_{i\in\bZ/m\bZ} \fg_i
\eeqn
where $\Lg_i$ is the $\zeta_m^i$-eigenspace of $d\theta$ for a fixed primitive $m$-th root of unity $\zeta_m$.  In~\cite{VX2} we focused on (GIT) stable gradings, i.e., on gradings such the action of $K$ on $\Lg_1$ is stable in the sense of geometric invariant theory. Below we work with arbitrary gradings. 

 We recall some terminology introduced in \autoref{invariant} applied to the grade Lie algebras setting. Recall that the $K$-orbits of semisimple  elements $x\in \fg_1$ are closed. We say that an element $x\in \fg_1$ is {\em regular} if $\on{dim}Z_K(x)\leq \on{dim}Z_K(y)$ for all $y\in\Lg_1$.  We say that a semisimple element $x_s\in\Lg_1$ is {\em regular semisimple} if $\on{dim}Z_K(x_s)\leq \on{dim}Z_K(y_s)$ for all semisimple elements $y_s\in\Lg_1$. It is important to note that it might happen that a regular semisimple element is not regular, but is merely regular among semisimple elements. We will use the following notation. We write $\fg_1^{\mathrm{s}}$ for the set of semisimple elements, $\fg_1^{\mathrm{rs}}$ for the set of regular semisimple elements, and $\fg_1^{\mathrm{r}}$ for the set of regular elements. Finally, we write $\fg_1^{\reg}$ for the set of strongly regular elements, that is, regular elements whose semisimple part is  regular semisimple (as in \autoref{invariant}).  
 
Let $\Char_K(\Lg_1)$ denote the set of  character sheaves on $\Lg_1$, that is, irreducible $K$-equivariant $\bC^*$-conic perverse sheaves $\cF$ on $\Lg_1$ whose singular support is nilpotent. We write $\fg^{\on{nil}}$ for the nilpotent cone in $\fg$.  We identify $\Lg_1^*$ with $\Lg_{m-1}=\Lg_{-1}$ using a $K$-invariant and $\theta$-invariant non-degenerate symmetric bilinear form on $\Lg$. Then, by definition, the character sheaves are Fourier transforms of irreducible $K$-equivariant perverse sheaves on $\fg^{\on{nil}}_{-1}= \fg^{\on{nil}}\cap \fg_{-1}$.

We fix a Cartan subspace $\La$ of $\Lg_1$; i.e., $\La \subset \Lg_1$ is a maximal abelian subspace consisting of semisimple elements. By the rank of the graded Lie algebra we mean the dimension of the Cartan subspace $\fa$, and by the ``semisimple" rank we mean $\dim\fa-\dim\fa^K$. We write  $W_\La = W(K,\La)$ for the ``little" Weyl group:
\beqn
W_\La = W(K,\La) \coloneqq N_K(\La) / Z_K(\La).
\eeqn

Consider the group $M=Z_G(\fa)$. By definition it is $\theta$-stable  and hence gives us a graded Lie algebra in its own right. We have
\beqn
\fm=Z_\fg(\fa) =  Z_\fg(\fa)_{der} \oplus Z(Z_\fg(\fa))\,.
\eeqn
Note that $Z(Z_\fg(\fa)) = \fa \oplus \fa^\perp$ with $\fa^\perp$ a $\theta$-stable complement. We write $\Lu=Z_\fg(\fa)_{der}$ and then
\beqn
\fm_1 = Z_{\fg_1}(\fa) =  \Lu_1 \oplus \fa
\eeqn
and 
\beqn
\fg_1 =   Z_\fg(\fa)_1  \oplus [\fk,\fa] = \fa \oplus [\fk,\fa]  \oplus \Lu_1\,.
\eeqn
 Both decompositions  are $N_K(\fa)$-invariant. The factor $\Lu_1$ is a prehomogenous vector space for the $Z_K(\fa)^0$-action.  Note that $(K,\Lg_1)$ is stable as a polar representation if and only if $\Lu_1=0$.

\subsection{Equivariant geometry}  In this paper our goal is to construct cuspidal character sheaves, i.e., those that do not arise as a direct summand (up to shift) from parabolic induction of character sheaves on graded $\theta$-stable Levi subalgebras contained in proper $\theta$-stable parabolic subgroups.  We denote the set of cuspidal character sheaves by $\Char_K^{\on{cusp}}(\Lg_1)$. In the cases we consider here, all cuspidal sheaves turn out to have full support, i.e., they are intersection homology extensions of irreducible $K$-equivariant local systems on $\fg_1^{\reg}$. We describe the equivariant fundament group $\pi_1^{K}(\fg_1^{\reg},v_0)$, $v_0\in\fg_1^{\reg}$ in this subsection.

We can associate to the pair $(G,\theta)$ a $\theta$-stable reductive subgroup $L\subset G$ such that $\Ll=\on{Lie}L$ contains $\fa$ in its derived algebra, the pair $(L^\theta=K_L,\Ll_1)$ is a stable polar representation,   and such that $(K,\Lg_1)$ and $(K_L,\Ll_1)$ have isomorphic rings of invariants. Thus  we have
\beq
 \fl_1 \inv K_L  \ = \ \fg_1 \inv K \ = \fa/W_\fa \ \,.
\eeq
We write  $$f_L:\fl_1 \to  \fa/W_\fa.$$
\begin{remark}
This fact does not appear to have a uniform proof but has been established by a case by case verification~\cite{L1,L2,RLYG}. In the cases we use it in this paper it can also be directly verified. 
\end{remark}

 \begin{remark}
 In \autoref{cuspidal sheaves}, where we apply the theory, the pairs $(K_L,\Ll_1)$ are stable in the  sense of graded Lie algebras, i.e., they are stable even in the GIT sense when considering the action of $K_L$ modulo its center. 
 \end{remark}

  Consider the adjoint quotient map
\beqn
f: \Lg_1 \to \Lg_1 \inv K \cong \La \slash W_\La \, .
\eeqn
It restricts to a fibration $$\tilde f: \Lg_1^{\reg} \to \La^{\rrs} \slash W_\La.$$ Let $v_0\in\fa^{\rrs}\times\Lu_1^\rr$,  and let $a_0$ be its component in $\fa^{\rrs}$ and $u_0$ its component in $\Lu_1^\rr$. Set $\barbasepta = f (v_0)$. Recall that
\beqn
B_{W_\La} \coloneqq \pi_1 (\La^{\rrs} \slash W_\La, \barbasepta)
\eeqn
is the braid group associated to the complex reflection group $W_\La$. We write 
\beq
 F^{\reg}=\tilde F_{\bar a_0}=  \tilde f^{-1} (\barbasepta) = K\cdot v_0\,
\eeq
for the fibre. It is an open subset in the fibre  $F =f^{-1} (\barbasepta)$. The fibre $F$ consists of a finite number of $K$-orbits with the closed orbit consisting of regular semisimple elements, i.e., the closed orbit is $F^{\rrs}=f^{-1} (\barbasepta)^{\rrs}$. We also note that  $f_L^{-1} (\barbasepta)\subset F^{\rrs}$. 

Recall that, following~~\cite{L1,L2,RLYG}, in~\cite{VX2} we construct the Kostant slice explicitly for $(K_L,\Ll_1)$ as follows. Consider a fundamental torus $T_{\mathbf{f}}$ of $L$, i.e., $T_{\mathbf{f}}$ is a $\theta$-stable maximal torus of $L$  and  $(T^\theta_{\mathbf{f}})^0$ is a maximal torus in $K_L$. Fix a pinning $(B,T_{\mathbf{f}},\{X_{\alpha}\}_{\alpha\in\Delta_L})$ of $L$ where we have written $\Delta_L$ for the simple roots of $(L,T_{\mathbf{f}})$. Consider the regular nilpotent element 
\beqn
E_1 = \sum_{\alpha \in \Delta_L} X_\alpha\in\fl_1\,.
\eeqn
Let  $\Ls_1 \subset \fl_1$ be  a slice containing $E_1$ such that the tangent space $[\fl_0,E_1]$ to the nilpotent cone of $\fl_1$ is transverse to $\Ls_1$ in $\fl_1$. We can now lift this to a Kostant slice for $(G,\theta)$  by choosing a generic point $E_2\in\Lu_1^\rr$ and then setting $E=E_1+E_2$. Proceeding as above consider a slice $\Ls\subset \Lg_1$ containing $E$ such that the tangent space $[\fg_0,E]$  to the nilpotent cone is transverse to $\Ls$ in $\fg_1$. As in~\cite{VX2} we conclude that the Kostant slice gives us a regular splitting in the terminology of~\cite{GVX2}. This implies that the exact sequence 
\beq
1 \to \pi_1^K( F^\reg,v_0 ) \to \pi_1^K(\Lg_1^{\reg},v_0) \to B_{W_\La} \to 1\,
\eeq
is split by the Kostant slice.  We use the same notation as in our previous papers and so we write $$\tilde B_{W_\La} =\pi_1^K(\Lg_1^{\reg},v_0).$$ Furthermore, we have an exact sequence
\beq
1 \to \pi_1^{Z_K(\fa)^0}( \Lu_1^\rr,u_0 ) \to \pi_1^K( F^{\reg},v_0 ) \to \pi_1^K( F^{\rrs},a_0) \to 1
\eeq
and we note that
\beq
\label{fiber fund group} 
 \pi_1^K( F^\reg,v_0 ) \cong \  \pi_1^{Z_K(\fa)}( \Lu_1^\rr,u_0 )\,.
\eeq

Writing $$I^{\rrs}=Z_K(\fa)/Z_K(\fa)^0= Z_{K_L}(\fa)/Z_{K_L}(\fa)^0$$  and $\tilde B_{W_\La}^{\rrs} =\pi_1^{K_L}(\Ll_1^{\rrs},a_0)$, 
we have, as in \cite[section 2.2]{GVX2}, a commutative diagram:
\beq
\label{lmainDiagram}
\begin{CD}
1 @>>> I^{\rrs} @>>> \widetilde B_{W_\La}^{\rrs} @>{\tilde q}>> B_{W_\La} @>>> 1 \;\,
\\
@. @| @VV{\tilde p}V @VV{p}V @.
\\
1 @>>> I^{\rrs} @>>> \widetilde W_\La^{\rrs}  @>{q}>> W_\La @>>> 1 \, 
\end{CD}
\eeq
where $ \widetilde W_\La^{\rrs} = N_K(\fa)/Z_K(\fa)^0$. In this diagram the top row is split by the Kostant slice and the action of $B_{W_\fa}$ on $I^{\rrs}$ factors through the obvious $W_\fa$-action.  
We will now write the corresponding diagram in the current setting where it is a bit more complicated.

First note that we have 
\beq
\tilde B_{W_\La} =  \pi^{N_K(\fa)}_1(\fa^{\rrs}\times\Lu_1^{\rr})\quad\text{ and }\quad   \tilde B_{W_\La}^{\rrs}= \pi^{N_K(\fa)}_1(\fa^{\rrs}) \,.  
\eeq
Consider the following commutative diagram with exact rows and with the right two columns also exact:
\beqn
\begin{CD}
@. I^0= \pi_1^{Z_K(\fa)^0}(\Lu_1^{\rr})@. \pi_1(\fa^{\rrs}) @= \pi_1(\fa^{\rrs})  @.
\\
@.@VVV   @VVV @VVV     @.
\\
1 @>>>I = \pi_1^{Z_K(\fa)}(\Lu_1^{\rr})@>>>\pi^{N_K(\fa)}_1(\fa^{\rrs}\times\Lu_1^\rr)  @>{\tilde q}>> \pi_1(\fa^{\rrs}/W_\fa)@>>> 1 \;\,
\\
@.@|   @VVV @VVV     @.
\\
1 @>>> I = \pi_1^{Z_K(\fa)}(\Lu_1^{\rr})@>>> \pi^{N_K(\fa)}_1(\Lu_1^\rr) @>{q}>> W_\fa @>>> 1 \, .
\end{CD}
\eeqn
Setting $\widetilde W_\La =  \pi^{N_K(\fa)}_1(\Lu_1^\rr)$ we obtain a diagram analogous to~\eqref{lmainDiagram}:
\beq
\label{mainDiagram}
\begin{CD}
1 @>>> I @>>> \widetilde B_{W_\La} @>{\tilde q}>> B_{W_\La} @>>> 1 \;\,
\\
@. @| @VV{\tilde p}V @VV{p}V @.
\\
1 @>>> I @>>> \widetilde W_\La @>{q}>> W_\La @>>> 1 
\end{CD}
\eeq
with the top row split by the Kostant slice. 

In our examples the group $I$ is abelian. Thus we obtain an action of $W_\fa$ on $I$, which is the natural action of $W_\La$ by conjugation on $I$. We have the following exact sequence of abelian groups with  a  $W_\fa$-action
\beq
\label{I's}
1 \to I^0 \to I \to I^\rrs \to 1\,.
\eeq
The action of $W_\fa$ on $I^0$ is trivial. As in the case of the existence of the stable $L$, we only establish this fact by case-by-case considerations. Furthermore, again arguing case-by-case, the exact sequence~\eqref{I's} can be naturally split. The verification of these facts will be done in the cases we consider in \autoref{cuspidal sheaves}.

The discussion above shows that the theory of character sheaves for general gradings can be split into a ``stable part" and an unstable ``semisimple rank zero" part. In particular, the Weyl group combinatorics in this paper is exactly the same as in our paper~\cite{VX2} on the stable gradings case. 

\begin{remark}
The hypothesis that $I$ is abelian is not essential. We always have an action of  $W_\fa$ on $\on{Rep}(I)$, the representations of $I$. For our constructions we only use this action. This generality is essential in the treatment of spin groups in~\cite{X}. We will also treat general gradings for spin groups in a future publication.
\end{remark}

\section{The topology of the nearby cycle construction}
\label{nearby}

In this section we analyse the nearby cycle construction which we briefly defined in \autoref{invariant}. 
Consider a graded Lie algebra satisfying the conditions in the last section and we continue to use the notation there. We make one further hypothesis:
\beq
\label{affine hypothesis}
\text{The stabilizer of $Z_K(\fa)^0$ of a generic point in $\Lu_1^\rr$ is reductive.}
\eeq
This condition implies that $\Lu_1$ is a finite regular prehomogenous vector space for $Z_K(\fa)^0$ and therefore $\Lu_1^\rr$ is affine.

 Consider a character $\chi$ of $I\cong Z_K(x)/Z_K(x)^0$, $x\in\Lg_1^{\reg}$. It gives rise to a rank one $K$-equivariant local system $\cL_\chi$ on $ F^{\reg}$.  We will consider its $\IC$-extension and we will assume that 
\beq\label{assumption clean}
\IC(F^\reg, \cL_\chi) \ \ \text{is clean}\,.
\eeq

Proceeding as in~\cite{GVX2} we can form the corresponding nearby cycle sheaf $P_\chi$. As explained in~\cite{GVX2} we take the limit of $\IC(\cL_\chi)$ to the nilpotent cone via the nearby cycle functor working on the dual space $\fg_{-1}$ and so we obtain a $K$-equivariant perverse sheaf $P_\chi$ on $\Lg^{\on{nil}}_{-1} = \fg_{-1} \cap \Lg^{\on{nil}}$.  We now consider the Fourier transform $\widehat P_\chi$. 

We recall the discussion in \autoref{invariant}. In our current notation $V= \fg_{-1}$, $V^*= \fg_{1}$, $(V^*)^{\reg}=\fg_1^{\reg}$, and the group $G$ in that section is now denoted by $K$. What is called in \autoref{invariant} the Cartan subspace $\fa$ is now $\fa_{-1}$ and the $\fa^*$ of \autoref{invariant} is now simply denoted by $\fa=\fa_1$.

Recall that $\widehat P_\chi  = \IC( \widehat P_\chi|_{\fg_1^{\reg}})$.
Let $\xi\in\fa^{\rrs} \times \Lu_1^\rr$. Then, after choosing appropriate cuts, we can, as explained in \autoref{invariant}, we express the stalk $(\widehat P_\chi)_\xi$ in the following terms:
\beq
(\widehat P_\chi)_\xi\ = \ \bigoplus_{w\in W_\fa} TM_w \otimes \mu_{F^{\rrs}}(\IC(\cL_\chi))_\xi\,.
\eeq
Here the $TM_w$ denotes the tangential Morse group at $w\cdot x_0^{\rrs}$. Using the Picard-Lefschetz cuts, as explained in~\cite{GVX1,GVX2}, the Morse groups represent elements in $(\widehat P_\chi)_\xi$. Note that when $\xi$ moves in $\fa^{\rrs} \times \Lu_1^\rr$, the first term in the tensor product in~\eqref{Morse data} does not change but the second term does      depend on  the  $\Lu_1^\rr$ component of $\xi$. The isomorphism in~\eqref{Morse data} depends on the  Picard-Lefschetz cuts which depend continuously on $\xi$.

In~\cite{GVX2} we made a careful analysis of the resulting local systems analysing explicitly the family of Hessians. The arguments in~\cite{GVX2} carry over verbatim to our setting and control the factor $TM_w$. In that paper we work in the stable polar situation in which case the entire fibre $F=F^{\rrs}$,  but that is the only distinction. There the factor $\mu_{F^{\rrs}}(\IC(\cL_\chi))_\xi$ is  a local system on $F=F^{\rrs}$. In~\cite{GVX2} we make certain hypotheses on the polar representations, i.e., visibility (or rank 1), locality, and the existence of regular nilpotents. The same arguments as in~\cite{VX2} show that these hypotheses are satisfied in the current setting as well. Thus we can reduce the description of $(\widehat P_\chi)_\xi$ to rank one.

The key in the explicit description of $ \widehat P_\chi|_{\fg_1^{\reg}}$ is reduction to rank one. Recall that the combinatorics in our setting is the same as for the polar stable pair $(L,\theta)$. As pointed out above this combinatorics describes the tangential part of the Morse data.

Let us choose a distinguished reflection $s\in W_\fa$ (such reflections are in bijection with reflection hyperplanes of $W_\fa$ in $\fa$). Let $\fa_s$ denote the corresponding reflecting hyperplane and let $G_s=Z_G(\fa_s)$. The automorphism $\theta$ preserves $G_s$ and so we have a new graded Lie algebra $\fg_s=\Lie(G_s)$ with $K_s = Z_K (\La_s)^0$. We then have
\beqn
\fg_{s,1} =   Z_\fg(\fa)_1  \oplus [\fk_s,\fa] = \fa \oplus [\fk_s,\fa]  \oplus \Lu_1\,.
\eeqn
This shows, in particular, that we get a semisimple rank one polar representation $(K_s,\fg_{s,1})$ with the same unstable part $\Lu_1$ as $(K,\Lg_1)$. We note that $G_s$  is not semisimple, but the results in~\cite[\S 3]{VX2} allow us to reduce to that case. It is now important to note that we have
\beq
\label{unstable does not change}
Z_K(\fa)^0 = Z_{K_s}(\fa)^0
\eeq
and so after reduction to rank one the unstable part given by the prehomogeneous vector space $(Z_K(\fa)^0 ,\Lu_1)$ does not change. 

We now proceed as in~~\cite{VX2}. So, we set
\beqn
W_{\fa,s} = N_{K_s} (\La) / Z_{K_s} (\La) \;\; \text{and} \;\; \widetilde W_{\fa,s} =  \pi^{N_{K_s}(\fa)}_1(\Lu_1^\rr).
\eeqn
Recall from~\cite{GVX2} that  $W_{\fa,s} = Z_{W_\fa}(\fa_s)=\langle s\rangle$, the subgroup of $W_\fa$ generated by the distinguished reflection $s$.
We have a natural projection $q_s : \widetilde W_{\fa,s} \to W_{\fa,s}$. Therefore, we have inclusions $W_{\fa,s} \to W_\La$ and $\widetilde W_{\fa,s} \to \widetilde W_\La$. We will use these inclusions to view $W_{\fa,s}$ (resp. $\widetilde W_{\fa,s}$) as a subgroup of $W_\La$ (resp. $\widetilde W_\La$). We have the following short exact sequence:
\beqn
1 \longrightarrow I_s \coloneqq  \pi_1^{Z_{K_s}(\fa)}(\Lu_1^\rr) \longrightarrow \widetilde W_{\fa,s} \xrightarrow{\; q_s \;} W_{\fa,s} \longrightarrow 1 \, .
\eeqn
Note that, by construction, we have $I_s  \subset I$.

Recall also that specifying a  splitting $\tilde r : B_{W_\La} \to \widetilde B_{W_\La}$ is equivalent to specifying a map $r: B_{W_\La}\to \widetilde W_\La$ such that $q\circ  r = p$. As our splitting comes from a Kostant slice it is regular and hence $r(\sigma_s)\in \widetilde W_{\fa,s}$, where $\sigma_s\in B_{W_\fa}$  is a chosen braid generator corresponding to $s$.

We set
\beqn
W_{\fa,\chi} \coloneqq \on{Stab}_{W_\fa} (\chi) \;\; \text{and} \;\;
B_{W_\fa}^\chi \coloneqq \on{Stab}_{B_{W_\fa}} (\chi) = p^{-1} (W_{\fa,\chi})\,.
\eeqn
For each distinguished reflection $s$, set 
\beqn
\text{$W_{\fa,s,\chi} = W_{\fa,s} \cap W_{\fa,\chi}$
and $e_s = |W_{\fa,s}|\ / |W_{\fa,s,\chi} | \in \bN$,}
\eeqn
 so that 
$W_{\fa,s,\chi} = \langle s^{e_s} \rangle$.  
Let $W_{\fa,\chi}^0 \subset W_{\fa,\chi}$ be the subgroup generated by all
the $W_{\fa,s,\chi}$. By definition $W_{\fa,\chi}^0$
is a complex reflection group acting on $\Cartan$.

Recall that in our setting the $W_\fa$-action on $I^0$ is trivial. As the exact sequence~\eqref{I's} is naturally split, we obtain characters $\chi^0$ of $I^0$ and $\chi^\rrs$ of $I^\rrs$ such that $\chi = \chi^0 \cdot \chi^\rrs$. Therefore,  $W_{\fa,\chi} =  W_{\fa,\chi^\rrs}$, and, as was pointed out earlier, the combinatorics in our current setting is the same as in the case of stable gradings. 

The character $\chi$ induces a character $\chi_s :I_s \to \bC^*$ for every distinguished reflection $s$, that is, $\chi_s=\chi|_{I_s}$. By~\eqref{unstable does not change} and the discussion below it we have $I_s^0=I^0$ and thus we obtain a $Z_{K_s}(\fa)^0$-equivariant local system on $\Lu_1^r$. We note again that this local system is the same local system we had before passing to rank one. 

We can now proceed as in~\cite{GVX2,VX2}.  To each character $\chi$ as above we  associate a polynomial 
$$
R_{\chi,s}\in\bC[x],
$$
which is the minimal polynomial of the monodromy around the hyperplane $\fa_s$.
As in~\cite{GVX2}  there exists $\bar R_{\chi,s}\in\bC[x]$ such that 
\beq\label{mono-1}
R_{\chi,s}(x)=\bar R_{\chi,s}(x^{e_s}),
\eeq
i.e. that $R_{\chi,s}$ is in fact a polynomial in $x^{e_s}$. 
 Note that the distinguished reflections in $W_{\fa,\chi}^0$ are  precisely the $s^{e_s}$. Thus, we can define a Hecke algebra $\cH_{W_{\fa,\chi}^0}$ associated to the complex reflection group $W_{\fa,\chi}^0$ with relations given by the polynomials $\bar R_{\chi,s}$.

Let
\beqn
B_{W_\fa}^{\chi, 0}=p^{-1}(W_{\fa,\chi}^{0}),\ \widetilde B_{W_\fa}^\chi = \tilde q^{-1} (B_{W_\fa}^\chi)
\;\; \text{and} \;\;
\widetilde B_{W_\fa}^{\chi, 0} = \tilde q^{-1} (B_{W_\fa}^{\chi, 0}).
\eeqn
As in~\cite{GVX2}, $ \widehat P_\chi|_{\fg_1^{\reg}}$ is the $K$-equivariant local system on $\Lg_1$ corresponding to the following representation of $\pi_1^{K}(\Lg_1^{\reg})=\widetilde{B}_{W_\fa}$
\beq
\label{nearby input}
 \widehat P_\chi|_{\fg_1^{\reg}}\ = \ \left( \bC [\widetilde B_{W_\fa}]
\otimes_{\bC [\widetilde B_{W_\fa}^{\chi, 0}]}
(\cL_{\chi} \otimes  \cH_{W_{\fa,\chi}^0}) \right)\otimes \bC_\tau \, ,
\eeq
where the Hecke algebra $\cH_{W_{\fa,\chi}^0}$ is viewed as a
$\bC [\widetilde B_{W_\fa}^{\chi, 0}]$-module via the composition of maps
\beqn
\widetilde B_{W_\fa}^{\chi, 0}\to B_{W_\fa}^{\chi, 0} \to
B_{W_{\fa,\chi}^0} \, ,
\eeqn
$\bC_\chi$ is viewed as a
$\bC [\widetilde B_{W_\fa}^{\chi, 0}]$-module via the unique extension of the character $\chi$ to $\widetilde B_{W_\fa}^{\chi, 0}$, and $\tau$ is the character on $I$ given by $\tau(x)= \det(x|_{\fg_1})$. As observed in~\cite{GVX2}, the character $\tau$ takes values in $\{\pm 1\}$.

Formula~\eqref{nearby input} reduces the description of $\widehat P_\chi|_{\fg_1^{\reg}}$ to the determination of the Hecke relations $\bar R_{\chi,s}$ in the Hecke algebra $\cH_{W_{\fa,\chi}^0}$. Furthermore, the calculation of the polynomials $\bar R_{\chi,s}$ is reduced to a  calculation in a rank one situation. 
\begin{remark}
\label{no tau}
In~\cite{GVX2} the formula~\eqref{nearby input} involves another character $\rho$ but the rank one calculations show that it is trivial in the setting of \autoref{cuspidal sheaves}. 
\end{remark}

\section{The Exact case}
\label{Exact}

In the rest of the paper we apply the considerations in the previous section but strengthen the condition~\eqref{affine hypothesis} as follows
\beq
\label{exact hypothesis}
\text{The stabiliser of $Z_K(\fa)^0$ of a generic point in $\Lu_1^\rr$ is finite.}
\eeq
 This is exactly the set-up for the existence of full support cuspidal character sheaves, as can be directly seen from~\cite[Theorem 3.7(i)]{LTVX}. Following~\cite{U}, such prehomogeneous vector spaces $(Z_K(\fa)^0,\Lu_1)$ are called exact. We will say that a grading on $\Lg$ is an exact grading if the condition~\eqref{exact hypothesis} holds.
 
 We view condition~\eqref{exact hypothesis} as a kind of stability condition. For GIT stability one assumes that generic obits are closed and generic stabilisers are finite. For stability for polar representation we drop the finite stabiliser condition but assume that generic orbits are closed. Our stability condition drops the requirement that generic orbits are closed but we retain the condition that generic stabilisers are finite. 

The task in \autoref{cuspidal sheaves} is to explicitly write down the list of full support character sheaves which will turn out to exhaust cuspidal character sheaves in the setting of classical groups. For this we need to explicitly determine the Hecke algebras  $\cH_{W_{\fa,\chi}^0}$. To do so we need to determine explicitly the polynomials $R_{\chi,s}$ obtained by reduction to rank one. We make use of b-functions as was explained in \autoref{sec-rank1g}. To facilitate this process we will spell out the translation of the results on invariant systems to the current graded Lie algebras setting. 

\subsection{The determinant character of an exact prehomogeneous vector space}

We use the notation of  \autoref{zero}. Let $(G,V)$ be a regular prehomogeneous vector space.  In this subsection we assume further that $\dim G = \dim V=n$, i.e., $(G,V)$ is an exact prehomogenous vector space. This is the analogue of hypothesis~\eqref{exact hypothesis}.

There always exists a semi-invariant corresponding to the character $(\on{det}\rho)^2$, see, for example, \cite[Corollary 2.17]{K}. However, for exact prehomogenous vector spaces we can do better. 
 In this case there exists a semi-invariant $f$ such that the corresponding character is $\on{det}\rho$ (see, for example, \cite[Corollary 2.20]{K}). This can be easily seen as follows. Pick $v\in V$ and consider the map $r_v: G \to V$ given by $g\mapsto gv$. Then we have $dr_v:\fg \to V$ and  $\wedge^n  dr_v \in \wedge^n V \otimes (\wedge^n \fg)^{-1}$. Identifying $\wedge^n V $ and $\wedge^n \fg $, the $\wedge^n  dr_v$ can be viewed as a function $f:V\to \bC$. By definition of $f$ we have 
\beq
f(gv) = \det (\rho(g)) f(v)\,,
\eeq
i.e., $f$ is a semi-invariant associated to the character $\det(\rho)$.

Recall that $f_i$ are the fundamental semi-invariants and $\chi_i$ are the corresponding characters of $G$. We claim that
\beq
\label{exact determinant}
f = \prod f_i \qquad \text{up to a constant and hence} \qquad \det(\rho)= \prod \chi_i\,.
\eeq
To see this we consider a general $v$ in the open orbit in the locus given by $f_i=0$. Because the orbit is of codimension one, the connected component of the stabilizer $G^v$ is $\bC^*$. By the Luna slice theorem~\cite{Lun} there is a one dimensional slice $S$ through $v$ (\'etale or complex analytic) such that on the slice $S$ locally near $v$ we can identify the action of $G^v$ on $S$ with a linear action of  $\bC^*$ on  $\bC$. It is easy to check that, up to a constant, the function $f$ coincides with the coordinate function $z$ on $\bC$. Hence the function $f$ vanishes to first order at the generic point of $f_i=0$. Thus, as both $f$ and $f_i$ are semi-invariants, we conclude that generically along $f_i=0$ the functions $f$ and $f_i$ differ by a constant. 

In particular we conclude:
\beq
\label{trivial tau}
\text{in the exact case the character $\tau$ in formula~\eqref{nearby input} is trivial.}
\eeq

\subsection{Rank one exact case} Let us consider the exact rank one situation, that is, the polar representation $(G,V)$ has rank one and  the generic stabiliser  in $\bC^*\times G$ is finite. We will formulate the conclusion of Theorem~\ref{thm-ft} in this case. By~\eqref{exact determinant} we see that  $\cL_{\chi\det(\rho)} =\cL_{\chi}$ and $\kappa_0$ is an integer. Thus the equations of Theorem~\ref{thm-ft} reduce to
\beq
\label{exact equations}
\widehat \cM_\chi|_{(V^*)^{\reg}} \ = \ \cL_{\chi} \otimes \bC[x]/ R_\chi(x) \qquad R_\chi(x) = \prod (x-\exp(2\pi i\alpha_i))^{d_i}\,.
\eeq

\subsection{The IC property} 
\label{Proof of IC} Recall that we make crucial use~\eqref{IC} in our proof in~\autoref{sec-rank1g}.  As the paper~\cite{GVX3} establishing this fact is not yet available we will briefly explain how to conclude~\eqref{IC} under the hypothesis~\eqref{exact hypothesis}. This is the context in which it is applied in the next section. The exactness condition~\eqref{exact hypothesis} implies that, in the language of~\cite{LTVX},  the locus $\fg_1^\reg$ is $G_0$-distinguished. The adjunction formula~\cite[Theorem 3.9]{LTVX} then implies that
\beqn
\dim\on{Hom} (\widehat P_\chi, \cG) =  \dim\on{Hom} (\bD\cG, \tau_*\IC(\cL_\chi)[\dim \fg_1 - \dim \fk])\,.
\eeqn
In this formula  $\cG$ is a character sheaf and $\tau$ is an identification of $\fg_{-1}$ and $\fg_1$ as explained in~\cite[section 3]{LTVX}. Now, as $\IC(\cL_\chi)$ is zero outside of $\fg_1^\reg$ the righthand side can only be non-zero if the support of $\cG$ is all of $\fg_1$. This means that $\widehat P_\chi$ cannot have any quotients whose support is not all of $\fg_1$. This implies that $\widehat P_\chi$ is an $\IC$-sheaf.

\section{Cuspidal Character Sheaves}
\label{cuspidal sheaves}

In this section we apply the theory developed in the previous sections to graded Lie algebras, more precisely, to Vinberg's type I gradings (see~\cite[\S7.2]{Vin}) on classical Lie algebras.  The inner type I gradings for special linear groups are dealt with in~\cite{X2}, see Conjecture 4.8 of {\em loc. cit.} for a conjectural description of all cuspidal character sheaves. In the remaining cases we  write down explicitly the  gradings which afford cupsidal character sheaves.  For each such gradings we write down a family of full support character sheaves, which are precisely the cuspidal character sheaves.  Furthermore, we describe these character sheaves explicitly in terms of representations of Hecke algebras. The fact that these sheaves are cuspidal follows directly from the main theorem in~\cite{LTVX}. The fact that they constitute all cuspidal character sheaves requires a minor argument beyond~\cite{LTVX}. We will address this issue in a future publication.

\subsection{Complex reflection groups and their Hecke algebras} The complex reflection groups and their Hecke algebras which will arise in  this section have been been discussed in~\cite[Subsection 7.1]{VX2}. In what follows we use the notation introduced there. The Weyl groups that  occur below will all be complex reflection groups of the form $G_{m,1,k}= (\bZ/m\bZ)^k \rtimes S_k$. The group $G_{m,2,k}$ came up in our work on the stable gradings case in~\cite{VX2}, but does not appear if we exclude the stable gradings. As in~\cite{VX2} we write  
$\cH^{a}G_{m,1,k}$ for the Hecke algebra of $G_{m,1,k}$ where the Hecke relations are $(x-1)^2=0$ for reflections of order 2 and  $(x-1)^a(x+1)^{m-a}=0$ for reflections of order $m$. 
 
\subsection{Vinberg's type I graded Lie algebras }For a positive integer $k$, we write $\zeta_k=e^{\frac{2\pi\mathbf{i}}{k}}$, a primitive $k$-th root of unity. 
Let $\theta:G\to G$ be an automorphism such that $\theta$ induces a type I grading of the Lie algebra $\Lg=\on{Lie}G$. Let $m$ be the order of $\theta$.  
We have  $$\Lg=\oplus_{k\in\bZ/m\bZ}\Lg_k,$$ where $\Lg_k$ is the eigenspace of $d\theta$ on $\Lg$ with eigenvalue $\zeta_m^k$.  
 We make use of the following concrete description of $\theta$ following Vinberg \cite{Vin} and Yun~\cite{Y}. We will omit the inner automorphisms of type I for special linear groups below as these cases were dealt with in~\cite{X2}. We will also write a superscript $^2$ in type A to indicate that we are considering the outer automorphisms of type I. However we will not seperate the inner and outer automorphisms for type D Lie algebras. 
 
{\bf Type  ${}^2$AI}  Let $G=SL_V$, where $V$ is a $\bC$-vector space of dimension $N$ equipped with a non-degenerate bilinear form $(-,-)$ (not necessarily symmetric). Let $$\theta:G\to G,\ \theta(g)=(g^*)^{-1},$$
where $g^*:V\to V$ is defined by $(gv,w)=(v,g^*w)$, for any $v,w\in V$. We have 
\beqn
\text{$\theta^2(g)=\gamma g\gamma^{-1}$ for some $\gamma\in G$.}
\eeqn
  Then $m=2d$ for some odd integer $d\in\bZ_{+}$.  We can and will choose $\gamma$ such that $\gamma^{d}= 1$.

 {\bf Type CI (resp. BDI)} Let $G=Sp_V$ (resp. $G=SO_V$), where $V$ is a $\bC$-vector space of dimension $2n$ (resp. $N$) equipped with a non-degenerate  symplectic (resp. symmetric bilinear) form $(\ ,)$.  Let
 \beqn
 \text{ $\theta:G\to G$, $\theta(g)=\gamma g\gamma^{-1}$, $\gamma\in Sp_V$ (resp. $O_V$).}
  \eeqn Then $m$ is even and we can and will assume that $\gamma^{m}=-1$ (resp. $\gamma^{m}=1$). 
 
 \subsubsection*{Quiver description}
Let $$V_\lambda=\{v\in V\mid \gamma v=\lambda v\},\,\lambda\in\bC.$$ Then $V=\bigoplus V_\lambda$. Note that $(V_\lambda,V_\mu)=0$ unless $\lambda\mu=1$. We have
 \begin{align*}
&\Lg_k=\{x\in\mathfrak{sl}_V\mid xV_\lambda\subset V_{\zeta_{d}^k\lambda},\ (xv,w)+\zeta_{m}^k(v,xw)=0,\ \forall v,w\in V\}&&\text{ type {$^2$}A}\\
&\Lg_k=\{x\in\mathfrak{sp}_V\text{ (resp. $\mathfrak{so}_V$)}\mid xV_\lambda\subset V_{\zeta_{m}^k\lambda}\}&&\text{ type C (resp. BD)}.
\end{align*}
Let 
\begin{align}\label{def-mi}
&l=(d-1)/2,&&M_i=V_{\zeta_{d}^{-i}},\ i\in[0,d-1], &&\text{ type $^2$A}\nonumber\\
 &l=m/2,&&  M_i=V_{\zeta_{2m}^{1-2i}},\ i\in[1,m], &&\text{ type C}\\
 &l=m/2,&&  M_i=V_{\zeta_{m}^{-i}},\,i\in[0,m-1], &&\text{ type BD}\nonumber. 
 \end{align}
By~\cite{Vin} we have
$$r=\on{dim}\fa=\min\{\on{dim}M_i\}$$ and
\begin{align*}
& W_\fa=G_{d,1,r} &&\text{ type $^2$A},\\
 &W_\fa=G_{m,1,r}&&\text{ type BCD} \text{ except in type D when  $\on{dim}M_0=\on{dim} M_l=r$,}
\\ &W_\fa=G_{m,2,r} && \text{ type D when  $\on{dim}M_0=\on{dim} M_l=r$}.
\end{align*}
Note that  in type $^2$A, $(M_i,M_j)=0$ unless $i+j\equiv 0\mod d$ and $(\,,\,)|_{M_{0}}$ is a non-degenerate symmetric bilinear form; in type C, $(M_i,M_j)=0$ unless $i+j\equiv 1\mod m$; in type BD, $(M_i,M_j)=0$ unless $i+j\equiv 0\mod m$ and $(\,,\,)|_{M_{0}}$, $(\,,\,)|_{M_{l}}$ is non-degenerate.

We also note that in type $D$, the automorphism is inner if and only if both $\dim M_0$ and $\dim M_l$ are even.

 Following~\cite{Y}, we have the following description of $(K,\Lg_1)$ using cyclic quivers.

\noindent{\bf Type $^2$AI.} We can identify $\Lg_1$ with the set of representations of the  quiver \beqn
\xymatrix{&M_1\ar[dl]_{x_1}&M_2\ar[l]_{x_2}&\cdots\ar[l]&M_{l}\ar[l]_{x_{l}}\\M_0\ar[dr]_{x_1^*}&&&&&\\&M_{2l}\ar[r]_{x_2^*}&M_{2l-1}\ar[r]&\cdots\ar[r]_{x_{l}^*}&M_{l+1}\ar[uu]_{x_{l+1}}}
\eeqn
such that 
$
( x_i v_i,v_{d-i+1})+\zeta_m(v_i,x_i^*v_{d-i+1})=0,\,i\in[1,l],\,(x_{l+1}v_{l+1},w_{l+1})=( x_{l+1}w_{l+1},v_{l+1}),$ $\forall v_i,w_i\in M_i. 
$
We have 
\bern
&&K\cong SO_{M_0}\times \prod_{i=1}^lGL_{M_i}, \,\Lg_1\cong\bigoplus_{i=1}^l\on{Hom}(M_i,M_{i-1})\oplus\on{Sym}^2(M_l).
\eern
Let $g=(g_0,g_1,\ldots,g_l)\in K$,  and $x=(x_1,\ldots,x_l,x_{l+1})\in\Lg_1$, where $^tg_0g_0=1$, $g_i\in GL_{M_i}$, $x_i\in\on{Hom}(M_i,M_{i-1})$, $i\in[1,l]$, and $x_{l+1}\in\on{Sym}^2(M_l)$.  The action of $g$ on $x$ can be identified as follows
\beqn
x_i\mapsto g_{i-1}x_ig_i^{-1},\, \,1\leq i\leq l,\,x_{l+1}\mapsto g_lx_{l+1}{}^tg_l.
\eeqn
\textbf{Type CI.} We can identify $\Lg_1$ with the set of representations of the quiver
\beqn
\xymatrix{&M_1\ar[d]_{x_1}&M_2\ar[l]_{x_2}&\cdots\ar[l]&M_l\ar[l]_{x_l}\\&M_{2l}\ar[r]_{x_2^*}&M_{2l-1}\ar[r]&\cdots\ar[r]_{x_l^*}&M_{l+1}\ar[u]_{x_{l+1}}}
\eeqn
such that
$( x_i v_i,v_{m-i+2})+( v_i,x_i^*v_{m-i+2})=0,\,i\in[2,l]$,
$ ( x_1 v_1,w_1)= ( x_1w_1,v_1),\ 
( x_{l+1}v_{l+1},w_{l+1})=( x_{l+1}w_{l+1},v_{l+1}),\,\forall v_i,\,w_i\in M_i. 
$
We have
\bern
K\cong\prod_{i=1}^l GL_{M_i},\ \Lg_1\cong\on{Sym}^2(M_1^*)\oplus\bigoplus_{i=2}^l\on{Hom}(M_i,M_{i-1})\oplus  \on{Sym}^2(M_l).
\eern
Let $g=(g_1,\ldots,g_l)\in K$ and $x=(x_1,x_{2},\ldots,x_l,x_{l+1})\in\Lg_1$, where $g_i\in GL_{M_i}$, $i\in[1,l]$, $x_1\in\on{Sym}^2(M_1^*)$, $x_i\in\on{Hom}(M_i,M_{i-1})$, $i\in[2,l]$, $x_{l+1}\in\on{Sym}^2(M_l)$. The action of $g$ on $x$ is then given as follows
\beqn
x_1\mapsto (^tg_1)^{-1}x_1g_1^{-1},\,x_i\mapsto g_{i-1}x_ig_i^{-1},\,x_{l+1}\mapsto g_lx_{l+1}\,{}^tg_l.
\eeqn
\textbf{Type BDI.} We can identify $\Lg_1$ with the set of representations of the quiver
\beqn
\xymatrix{&M_1\ar[dl]_{x_1}&M_2\ar[l]_{x_2}&\cdots\ar[l]&M_{l-1}\ar[l]_{x_{l-1}}\\M_0\ar[dr]_{x_1^*}&&&&&M_l\ar[ul]_-{x_l}\\&M_{2l-1}\ar[r]_{x_2^*}&M_{2l-2}\ar[r]&\cdots\ar[r]_{x_{l-1}^*}&M_{l+1}\ar[ur]_{x_l^*}}\,
\eeqn
 such that
$
(x_i v_i,v_{m-i+1})+(v_i,x_i^*v_{m-i+1})=0,\ i=1,\ldots,l,\ \forall\, v_i\in M_i.
$ 
We have
\bern
K\cong SO_{M_0}\times \prod_{i=1}^{l-1}GL_{M_i}\times  SO_{M_l},\ \ \Lg_1\cong\bigoplus_{i=1}^l\on{Hom}(M_i,M_{i-1}).
\eern
Let $g=(g_0,g_1,\ldots,g_l)\in K$ and $x=(x_1,x_{2},\ldots,x_l)\in\Lg_1$, where $^tg_0g_0={}^tg_lg_l=1$, $g_i\in GL_{M_i}$, $i\in[1,l-1]$, $x_i\in\on{Hom}(M_i,M_{i-1})$, $i\in[1,l]$. The action of $g$ on $x$ can be identified as $$
x_i\mapsto g_{i-1}x_ig_i^{-1},\,1\leq i\leq l.
$$

\subsection{Exact gradings and the automorphisms  \texorpdfstring{$\theta^{p,q}$}{theta p,q}}
\label{Exact gradings} In the rest of this section we will consider the following gradings which satisfy the condition~\eqref{exact hypothesis}. 

 Let $p,q$ be such that 
 \beq\label{eqn-pq}
 \text{$p,q\in\mathbb{Z}_{\geq 0}$, and $p+q\leq l$ in type $^2$A,\,B,\,D,\  $p+q\leq l-1$ in type C.}
 \eeq We consider the order $m$ automorphisms $\theta^{p,q}$ of  $G$ which induce gradings on $\Lg$ such that
 \beq\label{def of theta}
 \dim M_i=\begin{cases}q-i+r&0\leq i\leq q-1\qquad\quad\ \text{ type $^2$A,\,B,\,D}\\q-i+r+1&1\leq i\leq q\qquad\qquad\quad \text{ type C}\\r&q\leq i\leq l-p\quad\qquad\ \,\text{ type $^2$A,\,B,\,D}\\r&q+1\leq i\leq l-p\quad\ \ \text{ type C}\\r+i-l+p&l-p+1\leq i\leq l\end{cases}
 \eeq
where $r\in\bZ_{\geq 0}$. That is, the dimension vector is, in type ${}^2\text{A,\,BD}$, 
$$(\dim M_0,\dim M_1,\cdots,\dim M_l)=(r+q,r+q-1,\ldots,r+1,r,\ldots,r,r+1,\ldots,r+p)$$
 and in type C,
$$(\dim M_1,\cdots,\dim M_l)=(r+q,r+q-1,\ldots,r+1,r,\ldots,r,r+1,\ldots,r+p)\,.$$ 
We have $$r=\on{dim}\fa,\ \dim\Lg_1=\dim K+r$$ and so~\eqref{exact hypothesis} holds.  Furthermore, 
\begin{align*}
&N=dr+p(p+1)+q^2&&\text{ type $^2$A},
\\
&2n=mr+p(p+1)+{q(q+1)}&&\text{ type C},
 \\
 & N=mr+p^2+q^2 &&\text{ type BD}.
\end{align*} 

The grading induced by $\theta^{p,q}$ is (GIT) stable for type $^2$A if $p=0$ and $q\in\{0,1\}$, for type C if $p=q=0$, and for type BD if $p\in\{0,1\}$ and $q\in\{0,1\}$. These cases have been handled in~\cite{VX2}. Thus, in what follows, we sometimes exclude these cases. Note also that in~\cite{VX2} we work with spin groups but in this paper we restrict attention to special orthogonal groups. Because of these restrictions the complex reflection groups  $G_{m,2,k}$ do not appear in our description of cuspidal character sheaves. 

\begin{remark}We will deal with the spin groups and exceptional groups in a future publication.
\end{remark}

From now on we assume that we are in the situation of this subsection and  $\theta=\theta^{p,q}$ for some pair $(p,q)$ as in~\eqref{eqn-pq}. 
\subsection{The group \texorpdfstring{$I$ and its characters}{I and its characters}}
\label{ssec-i}
 Let $x\in\Lg_1^{\reg}$ and $I=Z_K(x)$. We have
\begin{align*}
&\text{Type $^2$A}&&I\cong(\bZ/2\bZ)^{r+p+q-1},&& \begin{array}{ll}I^{\rrs}\cong(\bZ/2\bZ)^{r} &\text{ if $q\geq 1$}\\  I^{\rrs}\cong(\bZ/2\bZ)^{r-1}&\text{ if $q=0$}\end{array} \\
&\text{Type C}&&I\cong(\bZ/2\bZ)^{r+p+q},&& I^{\rrs}\cong(\bZ/2\bZ)^{r},&&\\
&\text{Type BD}&&I\cong(\bZ/2\bZ)^{r+p+q-2}\text{ if $p+q\geq 1$},&&I^{\rrs}\cong(\bZ/2\bZ)^{r}\text{ if $p\geq 1$ and $q\geq 1$}, \\&&&I\cong(\bZ/2\bZ)^{r-1}\text{ if $p=q=0$}, &&I^{\rrs}\cong (\bZ/2\bZ)^{r-1} \text{ if }\min\{p,q\}=0.
\end{align*}
We choose a set of generators of $I$ as follows.

 There exists an $x$-stable orthogonal  decomposition 
 \beq 
 \label{decomposition}
 V=V^1\oplus\cdots\oplus V^r\oplus U^1\oplus\cdots\oplus U^q\oplus W^1\oplus\cdots\oplus W^p
 \eeq
 of the vector space $V$ such that
  \begin{align*}
&V^i=\oplus_{j=0}^{d-1} V^i_j &&U^i=U^i_{0}\oplus_{j=1}^{i-1}(U^i_{j}\oplus U^i_{d-j}) &&W^i=\oplus_{j=0}^{i-1}(W^i_{l+j}\oplus W^i_{l+1-j}) &&\text{ type $^2$A}
\\
&V^i=\oplus_{j=1}^{m} V^i_j &&U^i=\oplus_{j=1}^{i}(U^i_{j}\oplus U^i_{m+1-j}) &&W^i=\oplus_{j=0}^{i-1}(W^i_{l-j}\oplus W^i_{l+1+j}) &&\text{ type C}
\\
&V^i=\oplus_{j=0}^{m-1} V^i_j &&U^i=U^i_{0}\oplus_{j=1}^{i-1}(U^i_{j}\oplus U^i_{m-j}) &&W^i=W^i_l\oplus_{j=1}^{i-1}(W^i_{l-j}\oplus W^i_{l+j}) &&\text{ type BD},
 \end{align*}
 where $V^i_j,U^i_j,W^i_j\subset M_j$, $\on{dim}V^i_j=\on{dim}U^i_j=\on{dim}W^i_j=1$, and such that 
 \beqn
 \text{$x|_{V^i}$ is semisimple with eigenvalues $\zeta_{m_0}^ja_i$, $a_i\in\bC^*$, $j\in[0, m_0-1]$, $a_i\neq a_j$,}
 \eeqn
 where $m_0=d$ in type $^2$A and $m_0=m$ otherwise, 
  and 
  \begin{align*}
 &\text{$x|_{U^i}$ (resp. $x|_{W^i}$) is nilpotent with Jordan block of size $2i-1$ (resp. $2i$)}&&\text{type $^2$A}\\
 &\text{$x|_{U^i}$, $x|_{W^i}$ is nilpotent with Jordan block of size $2i$} &&\text{type C }\\
 &\text{$x|_{U^i}$, $x|_{W^i}$ is nilpotent with Jordan block of size $2i-1$}&&\text{type BD}.
 \end{align*}
 That is, we have chosen a Cartan subspace $\fa$ so that $\fa\subset\oplus_{i=1}^r\mathfrak{gl}(V^i)$ (cf.~\cite[\S7.1]{Vin}),  a Jordan decomposition $x=x_s+x_n$ of $x$ such that $x_s\in\fa^{rs}$,  and $x_n$ has Jordan blocks of sizes
 \begin{align*}
& (1+3+\cdots+2q-1)+(2+4+\cdots +2p)&&\text{type $^2$A}\\
&(2+4+\cdots+2p)+(2+4+\cdots+2q)&&\text{type C}\\
&(1+3+\cdots+2q-1)+(1+3+\cdots +2p-1)&&\text{type BD}.
 \end{align*}
  
 We define $\gamma_i:V\to V$, $i\in[1,r-1]$ such that  $\gamma_i|_{U^j,W^j}=1$ and 
 \beqn
\gamma_i|_{V^i,V^{i+1}}=-1,\,\gamma_i|_{V^j}=1,\,j\notin\{ i,i+1\}\,.
\eeqn 
We define $\beta_i:V\to V$,  such that  $\beta_i|_{V^j,U^j}=1$, and
\begin{align*}
&\beta_i|_{W^i}=-1,\,&&\beta_i|_{W^j}=1,\,j\neq i&&i\in[1, p]&&\text{ type $^2$A,\,C};\\
&\beta_i|_{W^1,W^i}=-1,\,&&\beta_i|_{W^j}=1,\,j\notin\{1, i\}&&i\in[1, p-1]&&\text{ type BD}.
\end{align*}
We define $\alpha_i:V\to V$, such that  $\alpha_i|_{V^j,W^j}=1$, and 
\begin{align*}
&\alpha_i|_{U^1,U^{i+1}}=-1 &&\alpha_i|_{U^j}=1,\,j\notin\{ 1,i+1\} &&i\in[1, q-1]&&\text{ type $^2$A,\,BD}\\
&\alpha_i|_{U^{i}}=-1 &&\alpha_i|_{U^j}=1,\,j\neq i && i\in[1, q]&&\text{ type C}.
\end{align*}
Suppose that 
\begin{center}$q\geq 1$ in type $^2$A, and $\min(p,q)\geq 1$ in type BD.
\end{center}  We further define $\gamma_r:V\to V$ as follows, 
\begin{align*}
&\gamma_r|_{V^r, U^1}=-1 &&\gamma_r|_{V^j}=1,\,j\neq r &&\gamma_r|_{U^j}=1,\,j\neq 1 &&\gamma_r|_{W^k}=1&&\text{ type $^2$A}.\\
&\gamma_r|_{V^r}=-1 &&\gamma_r|_{V^j}=1,\,j\neq r &&\gamma_r|_{U^k,W^k}=1 && \ &&\text{ type C}\\
&\gamma_r|_{V^r, U^1,W^1}=-1 &&\gamma_r|_{V^j}=1,\,j\neq r &&\gamma_r|_{U^k,W^k}=1,\,k\neq 1 && \ &&\text{ type BD}.
\end{align*}
 
One checks easily that we have the following explicit descriptions:
\beq
\text{$I=\langle\gamma_i,\alpha_i,\beta_i\rangle$, $I^0=\langle\alpha_i,\beta_i\rangle$ and $I^{\rrs}\cong\langle\gamma_i\rangle$. }
\eeq
This gives us the desired splitting 
$
I\cong I^0\times I^{rs}
$ mentioned in \autoref{graded}.

\subsection{The cleanness assumption}
\label{cleanness assumption}Let us define
\beqn
{}^0\widehat{I}=\{\chi\in\widehat{I}\mid\chi(\alpha_i)=(-1)^{i},\,\chi(\beta_i)=(-1)^{i}\}.
\eeqn
We claim that the set of characters in $\hat I$ satisfying assumption~\eqref{assumption clean} is exactly ${}^0\widehat{I}$, that is,
\beq\label{clean locs}
\on{IC}(F^\reg,\cL_\chi)\text{ is clean if and only if }\chi\in {}^0\widehat{I}.
\eeq
By the discussion in  \autoref{nearby}, it suffices to check the above claim in rank 1 which we will do in the following subsection.

We further define $\chi_k\in{}^0\widehat{I}$ by
\beq\label{eqn-chik-def}
\chi_k(\gamma_i)=1,\ i\neq k,\ \chi_k(\gamma_k)=-1.
\eeq

\subsection{The little Weyl group and its action on $I$}\label{ssec-weylaction}In~\autoref{graded} we introduced a $\theta$-stable reductive subgroup $L$ of $G$ such that the pair $(K_L,\fl_1)$ has the same invariant theory as $(K,\fg_1)$.  We give an explicit construction of the $L$ in the situations we consider. In terms of the  subspaces  in~\eqref{decomposition}, it is the subgroup of $G$ that leaves the subspaces $\oplus_{i=1}^rV^i\oplus U^1_{0}$ (resp. $\oplus_{i=1}^rV^i$, $\oplus_{i=1}^rV^i\oplus U^1_0\oplus W^1_0$), $U^i_j$, $i\neq 1$, $W^i_j$ (resp. $U^i_j,W^i_j$; $U^i_j,i\neq 1,W^i_j,i\neq 1$) invariant in type $^2$A (resp. C, BD). Then $L$ is a $\theta$-stable reductive subgroup that satisfies the properties described in  \autoref{graded}. 
The resulting derived group is
\begin{align*}
&L_{der}=SL(\oplus_{i=1}^rV^i\oplus U^1_{0}) &&\text{ type $^2$A}
\\
&L_{der}\cong Sp(\oplus_{i=1}^rV^i) &&\text{ type C}
\\
&L_{der}\cong SO(\oplus_{i=1}^rV^i\oplus U^1_0\oplus W^1_0) &&\text{ type BD}.
\end{align*}
Here and above we treat  $U^1_{0}$ and $W^1_0$ as the zero vector space when they are not defined. 
We note that the resulting  pair $(K_L,\fl_1)$ is GIT stable, modulo centre, $\fl_1$ contains the chose Cartan subspace $\fa$, and the Weyl group $W_\fa$ is the same for the pairs $(K_L,\fl_1)$ and $(K,\fg_1)$. That is
$$W_\fa\cong N_{K_L}(\fa)/Z_{K_L}(\fa).$$Moreover, the group $I^{\rrs}\cong Z_{K_L}(\fa)/Z_{K_L}(\fa)^0$.

It is easy to see that $W_\fa$ acts trivially on $I^0$ by observing that for any $g\in N_{K_L}(\fa)$, $g$ leaves $U^1_0$ and $W^1_0$ invariant. Moreover, as the Weyl group $W_\fa$ acts trivially on $I^0$, in view of the splitting $I\cong I^0\times I^{rs}$, we can identify the action of $W_\fa$ on ${}^0\widehat{I}$ with the action of the Weyl group for  $(K_L,\fl_1)$ on $\widehat{I^{\rrs}}$. The latter action has been studied in detail in~\cite{VX2}.

\subsection{Rank 1}\label{ssec-rank1}
  Let $r=1$ in~\eqref{def of theta}. We are in the situation of exact rank 1 gradings. In this subsection we apply the discussion in \autoref{sec-rank1g} and \autoref{Exact} to the polar representation $(K,\Lg_1)$. We continue to use the notations in the previous sections.

By the discussion in the end of last subsection, the braid group $B_{W_\fa}\cong\mathbb{Z}$ acts trivially on $I=Z_K(x)$, $x\in \Lg_1^{\reg}$.  Recall from \autoref{sec-rank1g} that, in this situation, if $\chi_0,\chi_1,\ldots,\chi_{\ke}$ are the characters of $K$ corresponding to fundamental semi-invariants $f_i\in\on{SI}(K,\Lg_1)$, and $f_0$ occurs in the invariant $f\in\bC[\Lg_1]^K$ as a factor of multiplicity 1, then we have the following:
  \begin{lemma}
  We have a natural bijection
  \beqn
  \left\{\chi=\prod_{i=1}^{\ke}\chi_i^{s_i}\in X^*(K)\mid s_i\in\bQ/\bZ\right\}\xrightarrow{\sim} \widehat I,\ \ \chi\mapsto\chi|_I.
  \eeqn
  \end{lemma}
  We will make use of this lemma in this section. In what follows we use $\upsilon_i$ to denote the characters corresponding to fundamental semi-invariants (to avoid notation conflict with~\eqref{eqn-chik-def}).
\subsubsection{Type   {\rm $^2$AI}} The (fundamental) semi-invariants in $\on{SI}(K,\Lg_1)$ are given as follows:
\begin{align*}
&f_i(x)=\on{det}(^t(x_1\cdots x_{i-1}x_i)x_1\cdots x_{i-1}x_i),\,1\leq i\leq q  &&\upsilon_i(g)=\on{det}(g_i)^{-2}
\\
&f_i(x)=x_i,  \ \ q+1\leq i \leq l-p &&\upsilon_i(g)= g_{i-1} g_{i}^{-1}
\\
&f_i(x)=\on{det}(x_{i}x_{i+1}\cdots x_lx_{l+1}\,^t(x_{i}x_{i+1}\cdots x_l)),\,l-p+1\leq i\leq l &&\upsilon_i(g)= \on{det}(g_{i-1})^2
\\
&f_{l+1}=\on{det}(x_{l+1}) &&\upsilon_{l+1}(g)= \on{det}(g_{l})^2\,.
\end{align*}
We remark that $f_i$'s are irreducible unless $q=1$. In the latter case $f_1$ is a product of two linear factors.  The invariant is
\beqn
f=f_qf_{q+1}^2\cdots f_{l-p}^2f_{l-p+1}\text{ if $q\geq 1$ and $f=f_{1}^2\cdots f_{l-p}^2f_{l-p+1}$ if $q=0$}.
\eeqn

We have 
\beqn
\text{$I=\langle\gamma_1;\alpha_i,\,1\leq i\leq q-1;\beta_i,\,1\leq i\leq p\rangle$ when $q\geq 1$, and  $I=\langle\beta_i,\,1\leq i\leq p\rangle$ when $q=0$.}
\eeqn Recall that we have chosen $\fa$ such that $$\Lu_1=\{y\in\Lg_1\mid y{V^1}=0,\,yW^i_j\subset W^i_{j-1},\,yU^i_{j}\subset U^i_{j-1}\}.$$
Let $\tilde{f}_i=f_i|_{\Lu_1}$, $\tilde{f}^1=\prod_{i=1}^{q-1}\tilde{f}_i$, $\tilde{f}^2=\prod_{j=l-p+2}^{l+1}\tilde{f}_i$.
 By~\cite[Theorems I.1 and II.1]{Saf} we have
\bern
&&\partial_{\tilde{f}^1}\left((\tilde{f}^1)^s\prod_{i=1}^{q-1}\tilde{f}_i^{-s_i}\right)=2^{d_1}\tilde{b}^1(s)\left((\tilde{f}^1)^{s-1}\prod_{i=1}^{q-1}\tilde{f}_i^{-s_i}\right)\\
\label{bfn-u12}
&&\partial_{\tilde{f}^2}\left((\tilde{f}^2)^s\prod_{i=l-p+2}^{l+1}\tilde{f}_i^{-s_i}\right)=2^{d_2}\tilde{b}^2(s)\left((\tilde{f}^2)^{s-1}\prod_{i=l-p+2}^{l+1}\tilde{f}_i^{-s_i}\right)
\eern
where $d_i=\on{deg}\tilde{f}^i$, $i=1,2$, and 
\begin{subequations}
\beq\label{eqn-u1a1}
\tilde{b}^1(s)=\prod_{1\leq i\leq j\leq q-1}\prod_{b=1}^{j-i+1}\left((j-i+1)s-\sum_{a=i}^js_{q-a}+\frac{j-i+2}{2}-b\right)^2
\eeq
\beq\label{eqn-u1a2}
\begin{gathered}
\tilde{b}^2(s)=\prod_{1\leq i\leq j\leq p-1}\prod_{b=1}^{j-i+1}\left((j-i+1)s-\sum_{a=i}^js_{l-p+a+1}+\frac{j-i+2}{2}-b\right)^2\\
\hspace{.5in}\cdot \prod_{1\leq i \leq p}\prod_{b=1}^{p-i+1}\left((p-i+1)s-\sum_{a=i}^ps_{l-p+a+1}+\frac{p-i+2}{2}-b\right)\,.
\end{gathered} 
\eeq
\end{subequations}
Moreover, we have
\ber
\label{bfn-2}&&\partial_{f_q}(f_q^s\prod_{j=1}^{q-1}f_j^{s_j})=2^{2q}s^2\prod_{j=1}^{q-1}\left(s+\sum_{i=j}^{q-1}s_i+\frac{q-j}{2}\right)^2f_q^{s-1}\prod_{j=1}^{q-1}f_j^{s_j}\text{ (if $q\geq 1$)}\\
\label{bfn-1}&&\partial_{f_{l-p+1}}(f_{l-p+1}^s\prod_{j=l-p+2}^{l+1}f_j^{s_j})=2^{2p+1}s^2\left(s+\sum_{i=1}^ps_{l-p+1+i}+\frac{p}{2}\right)\\
&&\hspace{2in}\cdot\prod_{j=1}^{p-1}\left(s+\sum_{i=1}^js_{l-p+1+i}+\frac{j}{2}\right)^2f_{l-p+1}^{s-1}\prod_{j=l-p+2}^{l+1}f_j^{s_j}.\nonumber
\eer

Let us now choose $\chi_0=\upsilon_{l-p+1}$, here $\chi_0$ is as in~\autoref{sec-rank1g}. Let $$\chi=\prod_{j\in[1,l+1]\backslash\{l-p+1\}}\upsilon_j^{s_j}$$ and $u_\chi=\prod_{j\in[1,l+1]\backslash\{l-p+1\}}f_j^{-s_j}$. 
It follows that
\bern
&&\partial_f(f^su_\chi)=2^{d}b_{f,\chi}(s)f^{s-1}u_\chi,\eern
where $b_{f,\chi}(s)=b^1_{f,\chi}(s)b^2_{f,\chi}(s)b^3_{f,\chi}(s)$ with
\bern
\label{b1} &&b^1_{f,\chi}(s)=\prod_{j=1}^{q}\left(s-\sum_{a=j}^{q}s_a+\frac{q-j}{2}\right)^2\\
&&b^2_{f,\chi}(s)=s^2\left(s-\sum_{a=1}^ps_{l-p+1+a}+\frac{p}{2}\right)\prod_{j=1}^{p-1}\left(s-\sum_{a=1}^js_{l-p+1+a}+\frac{j}{2}\right)^2\nonumber\\&&
b^3_{f,\chi}(s)=\prod_{j=q+1}^{l-p}\left(\left(s-\frac{s_j}{2}\right)\left(s-\frac{s_{j}+1}{2}\right)\right).\nonumber
\eern
Here we use the convention that $b^1_{f,\chi}(s)=1$ if $q=0$, $b^2_{f,\chi}(s)=s$ if $p=0$, and $b^3_{f,\chi}(s)=1$ if $p+q=l$.
 
One readily checks that
\beqn
\chi\in X^*(K)\text{ if and only if } s_1,\ldots, s_q\in\frac{1}{2}\mathbb{Z},\ s_{q+1},\ldots,s_{l-p} \in\bZ,\text{ and }s_{l-p+2},\ldots,s_{l+1}\in\frac{1}{2}\mathbb{Z}\,.
\eeqn
Moreover, in view of~\eqref{eqn-u1a1},~\eqref{eqn-u1a2}, and Proposition~\ref{clean criterion}, we see that~\eqref{clean criterion rk1} holds if and only if 
\beqn
s_i\in\frac{1}{2}\mathbb{Z}\backslash\bZ,\,\text{ for } i\in[1, q-1]\cup[l-p+2, l+1].
\eeqn
Thus claim~\eqref{clean locs} follows, since we have
\begin{align*}
&\chi(\gamma_1)=(-1)^{-2\sum_{j=1}^qs_j+2\sum_{j=l-p+2}^{l+1}s_j},&& \text{when $q\geq 1$}, && \ &&&
\\
&\chi(\alpha_i)=(-1)^{-2\sum_{j=1}^is_j} &&1\leq i\leq q-1,&&& 
\\
&\chi(\beta_i)=(-1)^{2\sum_{j=l-i+2}^{l+1}s_j} &&1\leq i\leq p. &&&
\end{align*}

Let now 
$\chi=\left(\prod_{j=1}^{q-1}\upsilon_j^{1/2}\right)\upsilon_q^{s_q}\left(\prod_{j=l-p+2}^{l+1}\upsilon_j^{1/2}\right)$, where $s_q\in\{0,1/2\}$. Then $\chi \in{}^0\widehat I$.  For such a $\chi$ we have
\beqn
b_{f,\chi}(s)=s^{l+p-q+1}(s-\tfrac{1}{2})^{l-p-q}(s-s_q)^{2q}.\eeqn
Applying Theorem~\ref{thm-ft} in the form of~\eqref{exact equations} we conclude that
\beq
\label{rk1 AI poly}
\widehat{P_{\chi}}|_{\Lg_1^{\reg}}=\cL_\chi\otimes\bC[x]/(R_\chi(x)),\ \ R_\chi(x)=\begin{cases}(x-1)^{l+p+q+1}(x+1)^{l-p-q}\text{ if }s_q=0\\(x-1)^{l+p-q+1}(x+1)^{l-p+q}\text{ if }s_q=1/2\end{cases}
\eeq
where $\chi(\alpha_i)=(-1)^i ,\,i=1,\ldots,q-1,\ \chi(\beta_i)=(-1)^i,\,i=1,\ldots,p,$ and when $q\geq 1$, $\chi(\gamma_1)=(-1)^{p-q+1+2s_q}$ .

\subsubsection{Type {\rm CI}} The semi-invariants in $\on{SI}(K,\Lg_1)$ are given as follows:
\begin{align*}
&f_1(x)=\on{det}(x_1),&&\upsilon_1(g)=(\on{det}g_1)^{-2}\\
&f_i(x)=\on{det}({}^t(x_2\cdots x_{i-1}x_i)x_1x_2\cdots x_{i-1}x_i),\,2\leq i\leq q+1,&&\upsilon_i(g)=(\on{det}g_i)^{-2}\\
&f_i(x)=x_i,\,q+2\leq i\leq l-p,&&\upsilon_i(g)=g_{i-1}g_i^{-1}\\
&f_i(x)=\on{det}(x_{i}x_{i+1}\cdots x_l\,x_{l+1}\,^t(x_{i}\cdots x_l)),\,l-p+1\leq i\leq l,&&\upsilon_i(g)=(\on{det}g_{i-1})^{2}\\
& f_{l+1}(x)=\on{det}(x_{l+1}),&&\upsilon_{l+1}(g)=(\on{det}g_{l})^{2}.
\end{align*}
 The invariant is
\beqn
f=f_{q+1}f_{q+2}^2\cdots f_{l-p}^2f_{l-p+1}.
\eeqn
We have 
\beqn
\text{$I=\langle\gamma_1;\alpha_i,\,1\leq i\leq q;\beta_i,\,1\leq i\leq p\rangle$.}
\eeqn We have chosen $\fa$ such that $\Lu_1=\{y\in\Lg_1\mid y{V^1}=0,\,yW^i_j\subset W^i_{j-1},\,yU^i_{j}\subset U^i_{j-1}\}.$ 
Let $\tilde{f}_i=f_i|_{\Lu_1}$, $\tilde{f}^1=\prod_{i=1}^{q}\tilde{f}_i$, $\tilde{f}^2=\prod_{j=l-p+2}^{l+1}\tilde{f}_i$.
 By~\cite[Theorems I.1 and II.1]{Saf} we have
\bern
&&\partial_{\tilde{f}^1}\left((\tilde{f}^1)^s\prod_{i=1}^{q}\tilde{f}_i^{-s_i}\right)=2^{d_1}\tilde{b}^1(s)\left((\tilde{f}^1)^{s-1}\prod_{i=1}^{q}\tilde{f}_i^{-s_i}\right)\\
\label{bfn-u12}
&&\partial_{\tilde{f}^2}\left((\tilde{f}^2)^s\prod_{i=l-p+2}^{l+1}\tilde{f}_i^{-s_i}\right)=2^{d_2}\tilde{b}^2(s)\left((\tilde{f}^2)^{s-1}\prod_{i=l-p+2}^{l+1}\tilde{f}_i^{-s_i}\right)
\eern
where $d_i=\on{deg}\tilde{f}^i$, $i=1,2$, and 
\begin{subequations}
\beq\label{eqn-u1c1}
\begin{gathered}
\tilde{b}^1(s)=\prod_{1\leq i\leq j\leq q-1}\prod_{b=1}^{j-i+1}\left((j-i+1)s-\sum_{a=i}^js_{q+1-a}+\frac{j-i+2}{2}-b\right)^2\\\hspace{.5in}\cdot \prod_{1\leq i \leq q}\prod_{b=1}^{q-i+1}\left((q-i+1)s-\sum_{a=i}^qs_{q+1-a}+\frac{q-i+2}{2}-b\right)\,,
\end{gathered}\eeq
\beq\label{eqn-u1c2}
\begin{gathered}
\tilde{b}^2(s)=\prod_{1\leq i\leq j\leq p-1}\prod_{b=1}^{j-i+1}\left((j-i+1)s-\sum_{a=i}^js_{l-p+a+1}+\frac{j-i+2}{2}-b\right)^2\\
\hspace{.5in}\cdot \prod_{1\leq i \leq p}\prod_{b=1}^{p-i+1}\left((p-i+1)s-\sum_{a=i}^ps_{l-p+a+1}+\frac{p-i+2}{2}-b\right)\,.
\end{gathered}
\eeq
\end{subequations}

Let $\chi=\left(\prod_{j\in[1,l+1]\backslash\{q+1\}}\upsilon_j^{s_j}\right)$ (we have chosen $\chi_0=\upsilon_{q+1}$) and $u_\chi=\prod_{j\in[1,l+1]\backslash\{q+1\}}f_j^{-s_j}$.
It follows from~\eqref{bfn-1} that
\bern
&&\partial_f(f^su_\chi)=2^mb_{f,\chi}(s)f^{s-1}u_\chi,
\eern
where $b_{f,\chi}(s)=b^1_{f,\chi}(s)b^2_{f,\chi}(s)b^3_{f,\chi}(s)$ with
\ber
\label{b2}&& b^1_{f,\chi}(s)=s^{2}\left(s-\sum_{a=1}^qs_a+\frac{q}{2}\right)\prod_{j=2}^{q}\left(s-\sum_{a=j}^{q}s_a+\frac{q-j+1}{2}\right)^2\nonumber\\
&&b^2_{f,\chi}(s)=(s-s_{l-p+1})^2\left(s-\sum_{a=0}^ps_{l-p+1+a}+\frac{p}{2}\right)\prod_{j=1}^{p-1}\left(s-\sum_{a=0}^js_{l-p+1+a}+\frac{j}{2}\right)^2\nonumber\\&&b^3_{f,\chi}(s)=\prod_{j=q+2}^{l-p}\left(s-\frac{s_j}{2}\right)\left(s-\frac{s_j+1}{2}\right)\nonumber.
\eer
Here we use the convention that $b^1_{f,\chi}(s)=s$ if $q=0$, $b^2_{f,\chi}(s)=s-s_{l+1}$ if $p=0$, and $b^3_{f,\chi}(s)=1$ if $p+q=l-1$.

One readily checks that
\beqn
\text{$\chi\in X^*(K)$ if and only if }s_1,\ldots, s_q\in\frac{1}{2}\mathbb{Z},\,s_{q+2},\ldots, s_{l-p}\in\bZ,\text{ and }s_{l-p+1},\ldots, s_{l+1}\in\frac{1}{2}\bZ\,.
\eeqn
Moreover, in view of~\eqref{eqn-u1c1},~\eqref{eqn-u1c2}, and Proposition~\ref{clean criterion}, we see that~\eqref{clean criterion rk1} holds if and only if 
\beqn
s_i\in\frac{1}{2}\mathbb{Z}\backslash\bZ,\,\text{ for } i\in[1, q]\cup[l-p+2, l+1].
\eeqn
Thus claim~\eqref{clean locs} follows, since we have
\begin{align*}
&\chi(\gamma_1)=(-1)^{-2\sum_{j=1}^qs_j+2\sum_{j=l-p+1}^{l+1}s_j} && \ &&&
\\
&\chi(\alpha_i)=(-1)^{-2\sum_{j=1}^is_j} &&1\leq i\leq q,&&& 
\\
&\chi(\beta_i)=(-1)^{2\sum_{j=l-i+2}^{l+1}s_j} &&1\leq i\leq p. &&&
\end{align*}

Suppose now that  $\chi=\left(\prod_{j=1}^{q}\upsilon_j^{1/2}\right)\upsilon_{l-p+1}^{s_{l-p+1}}\left(\prod_{j=l-p+2}^{l+1}\upsilon_j^{1/2}\right)$, where $s_{l-p+1}\in\{0,1/2\}$. Then $\chi\in{}^0\widehat{I}$. For such a $\chi$ we have
\beqn
b_{f,\chi}(s)=s^{l-p+q}(s-\tfrac{1}{2})^{l-p-q-1}(s-s_{l-p+1})^{2p+1}.
\eeqn
Applying Theorem~\ref{thm-ft}  in the form of~\eqref{exact equations}, we conclude that
\beq
\label{rk1 CI poly}
\begin{gathered}
\widehat{P_{\chi}}|_{\Lg_1^{\reg}}=\cL_\chi\otimes\bC[x]/(R_\chi(x)),\\ R_\chi(x)=\begin{cases}(x-1)^{l+p+q+1}(x+1)^{l-p-q-1}\text{ if }s_{l-p+1}=0\\(x-1)^{l-p+q}(x+1)^{l+p-q}\ \ \ \ \ \ \text{ if }s_{l-p+1}=\frac{1}{2}\,.\end{cases}
\end{gathered}
\eeq
where 
$
\chi(\alpha_i)=(-1)^i ,\,i=1,\ldots,q,\ \chi(\beta_i)=(-1)^i,\,i=1,\ldots,p,\text{ and } \chi(\gamma_1)=(-1)^{p-q+2s_{l-p+1}}\, .
$

\subsubsection{Type {\rm BDI}}

The semi-invariants in $\on{SI}(K,\Lg_1)$ are given as follows:
\begin{align*}
&f_i(x)=\on{det}(^t(x_1\cdots x_{i-1}x_i)x_1\cdots x_{i-1}x_i),\,1\leq i\leq q&&\upsilon_i(g)=(\on{det}g_i)^{-2}\\
&f_i(x)=x_i,\,q+1\leq i\leq l-p&&\upsilon_i(g)=g_{i-1}g_i^{-1}\\
&f_i(x)=\on{det}(x_{i}x_{i+1}\cdots x_l\,^t(x_{i}x_{i+1}\cdots x_l)),\,l-p+1\leq i\leq l&&\upsilon_i(g)=(\on{det}g_{i-1})^{2}.
\end{align*}
As in type AI, the $f_i$'s are irreducible unless $q=1$ or $p=1$. In the latter cases we have $f_1$ or $f_l$ is a product of two linear factors, respectively. Let us write $f_0=f_{l+1}=1$. The invariant is
\beqn
f=f_qf_{q+1}^2\cdots f_{l-p}^2f_{l-p+1}.
\eeqn
Note that when $p,q\leq 1$, we are in the situation of stable gradings. These cases have been dealt with in~\cite{VX2}. 

Without loss of generality we assume that $q\geq 2$. We have 
\beqn
I=\langle\gamma_1 \text{ (when $p\geq 1$)};\,\alpha_i,\,1\leq i\leq q-1;\,\beta_i,\,1\leq i\leq p-1\rangle\,.
\eeqn We have chosen $\fa$ such that $\Lu_1=\{y\in\Lg_1\mid y{V^1}=0,\,yW^i_j\subset W^i_{j-1},\,yU^i_{j}\subset U^i_{j-1}\}.$
Let $\tilde{f}_i=f_i|_{\Lu_1}$, $\tilde{f}^1=\prod_{i=1}^{q-1}\tilde{f}_i$, $\tilde{f}^2=\prod_{j=l-p+2}^{l+1}\tilde{f}_i$.
 By~\cite[Theorems I.1 and II.1]{Saf} we have
\bern
&&\partial_{\tilde{f}^1}\left((\tilde{f}^1)^s\prod_{i=1}^{q-1}\tilde{f}_i^{-s_i}\right)=2^{d_1}\tilde{b}^1(s)\left((\tilde{f}^1)^{s-1}\prod_{i=1}^{q-1}\tilde{f}_i^{-s_i}\right)\\
\label{bfn-u12}
&&\partial_{\tilde{f}^2}\left((\tilde{f}^2)^s\prod_{i=l-p+2}^{l+1}\tilde{f}_i^{-s_i}\right)=2^{d_2}\tilde{b}^2(s)\left((\tilde{f}^2)^{s-1}\prod_{i=l-p+2}^{l+1}\tilde{f}_i^{-s_i}\right)
\eern
where
\begin{subequations}
\beq\label{eqn-u1bd1}
\tilde{b}^1(s)=\prod_{1\leq i\leq j\leq q-1}\prod_{b=1}^{j-i+1}\left((j-i+1)s-\sum_{a=i}^js_{q-a}+\frac{j-i+2}{2}-b\right)^2
\eeq
\beq\label{eqn-u1bd2}
\tilde{b}^2(s)=\prod_{1\leq i\leq j\leq p-1}\prod_{b=1}^{j-i+1}\left((j-i+1)s-\sum_{a=i}^js_{l-p+a+1}+\frac{j-i+2}{2}-b\right)^2\,.\eeq
\end{subequations}

 Let $\chi=\left(\prod_{j\in[1,l]\backslash\{q\}}\upsilon_j^{s_j}\right)$ (we have chosen $\chi_0=\upsilon_{q}$) and $u_\chi=\prod_{j\in[1,l]\backslash\{q\}}f_j^{-s_j}$.
It follows from~\eqref{bfn-2} that
\bern
&&\partial_f\left(f^su_\chi\right)=2^{m}b_{f,\chi}(s)f^{s-1}u_\chi\eern
where
\bern
b_{f,\chi}(s)&=&s^{2}\prod_{j=1}^{q-1}\left(s-\sum_{a=j}^{q-1}s_a+\frac{q-j}{2}\right)^2\prod_{j=0}^{p-1}\left(s-\sum_{a=0}^js_{l-p+1+a}+\frac{j}{2}\right)^2\\&&\prod_{j={q+1}}^{l-p}\left(s-\frac{s_j}{2}\right)\left(s-\frac{s_j+1}{2}\right).
\eern
One readily checks that  
\beqn
\text{$\chi\in X^*(K)$ if and only if }s_1,\ldots, s_{q-1}\in\frac{1}{2}\mathbb{Z},\,s_{q+1},\ldots, s_{l-p}\in\bZ,\text{ and }s_{l-p+1},\ldots, s_{l}\in\frac{1}{2}\bZ\,.
\eeqn
Moreover, in view of~\eqref{eqn-u1bd1},~\eqref{eqn-u1bd2}, and Proposition~\ref{clean criterion}, we see~\eqref{clean criterion rk1} holds if and only if 
\beqn
s_i\in\frac{1}{2}\mathbb{Z}\backslash\bZ,\,\text{ for } i\in[1, q-1]\cup[l-p+2, l].
\eeqn
Thus claim~\eqref{clean locs} follows, since we have
\begin{align*}
&\chi(\gamma_1)=(-1)^{-2\sum_{j=1}^qs_j+2\sum_{j=l-p+1}^{l}s_j} && \ &&&
\\
&\chi(\alpha_i)=(-1)^{-2\sum_{j=1}^is_j} &&1\leq i\leq q-1,&&&  \\
&\chi(\beta_i)=(-1)^{2\sum_{j=l-i+1}^{l}s_j} &&1\leq i\leq p-1. &&&
\end{align*}
Suppose now that $\chi=\left(\prod_{j=1}^{q-1}\chi_j^{1/2}\right)\chi_{l-p+1}^{s_{l-p+1}}\left(\prod_{j=l-p+2}^{l}\chi_j^{1/2}\right)$, where $s_{l-p+1}\in\{0,\tfrac{1}{2}\}$, if $p\geq 1$; and $\chi=\prod_{j=1}^{q-1}\chi_j^{1/2}$ if $p=0$. Then $\chi\in{}^0\widehat{I}$. For such a $\chi$ we have
\beqn
b_{f,\chi}(s)=s^{l-p+q}(s-\tfrac{1}{2})^{l-p-q}(s-s_{l-p+1})^{2p}. 
\eeqn
Applying Theorem~\ref{thm-ft}  in the form of~\eqref{exact equations}, we conclude that  
\beq
\label{rk1 BDI poly}
\widehat{P_{\chi}}|_{\Lg_1^{\reg}}=\cL_\chi\otimes\bC[x]/(R_\chi(x)),\ R_\chi(x)=\begin{cases}(x-1)^{l+p+q}(x+1)^{l-p-q}&\text{ if }s_{l-p+1}=0\\(x-1)^{l-p+q}(x+1)^{l+p-q} &\text{ if }s_{l-p+1}=\frac{1}{2}\,\end{cases}
\eeq
where
$
\chi(\beta_i)=(-1)^i,\,1\leq i\leq p-1,\,\chi(\alpha_i)=(-1)^i ,\,1\leq i\leq q-1,$  and  $\chi(\gamma_1)=(-1)^{p-q+2s_{l-p+1}} \text{ (if $p\geq 1$)}.
$

\subsection{Cuspidal character sheaves} In this section we write down explicitly the cuspidal character sheaves for the exact gradings in~\autoref{Exact gradings}. They  are IC-extensions of irreducible quotients of $\widehat P_\chi|_{\fg_1^\rrs}$ associated to the characters $\chi\in {}^0\widehat{I}$ in~\autoref{cleanness assumption} whose corresponding local system $\cL_\chi$ is clean.  To do that we will determine the local system $\widehat P_\chi|_{\fg_1^\reg}$ explicitly. To do so, we apply~\eqref{nearby input} and the reduction to rank 1 as explained in~\autoref{nearby}.

Recall from~\autoref{ssec-weylaction} that the Weyl group $W_\fa$ acts trivially on $I^0$ and that the groups $W_{\fa,\chi}$ for $(K,\fg_1)$ can be identified with analogous groups for  $(K_L,\fl_1)$. These groups are determined in ~\cite[Lemmas 7.12, 7.13 ]{VX2} for type {}$^2$A,~\cite[Lemma 7.15]{VX2} for type C,~\cite[Lemma 7.14, 7.17]{VX2} for type BD. 
Recall the characters $\chi_k\in{}^0\widehat{I}$ defined in~\eqref{eqn-chik-def}. Analysing the action of $W_\fa$ on ${}^0\widehat{I}$ and making use of the rank one calculations~\eqref{rk1 AI poly}, \eqref{rk1 CI poly} \eqref{rk1 BDI poly} we conclude that 
\begin{proposition} 
{\rm (i)} Suppose that $G=SL_N$. A set of representatives of $W_\fa$-orbits in ${}^0\widehat{I}$ is 
$$
\{\chi_k,\,0\leq k\leq r/2\}\text{ if }q=0,\ \ \{\chi_k,\,0\leq k\leq r\}\text{ if }q\geq 1.
$$
 We have
$$W_{\fa,\chi_k}^0\cong G_{d,1,k}\times G_{d,1,r-k},$$ and
\begin{align*}
  & W_{\fa,\chi_{r/2}}/ W_{\fa,\chi_{r/2}}^0\cong\bZ/2\bZ\,&&\text{ when $k= r/2$ and $q=0$}\\
  &W_{\fa,\chi_k}=W_{\fa,\chi_k}^0&&\text{ otherwise}.
\end{align*}
Moreover,  we have
\bern
&&\cH_{W_{\fa,\chi_k}^0}=\begin{cases}\cH^{l+p+q+1}(G_{d,1,k})\otimes \cH^{l+1+p-q}(G_{d,1,r-k})&\text{ if $p+q$ even}\\\cH^{l+1+p-q}(G_{d,1,k})\otimes\cH^{l+p+q+1}(G_{d,1,r-k}) &\text{ if $p+q$ odd}.\end{cases}
\eern

{\rm (ii)} Suppose that $G=Sp_{2n}$. A set of representatives of $W_\fa$-orbits in ${}^0\widehat I$ is $$\{\chi_k,\ 0\leq k\leq r\}.$$ We have that
$$
W_{\fa,\chi_k}=W_{\fa,\chi_k}^0\cong G_{m,1,k}\times G_{m,1,r-k},\ 0\leq k\leq r.
$$
Moreover, we have that
\bern
&&\cH_{W_{\fa,\chi_k}^0}=\begin{cases}\cH^{l-p+q}(G_{m,1,k})\otimes \cH^{l+p+q+1}(G_{m,1,r-k})&\text{ if $p+q$ even}\\\cH^{l+p+q+1}(G_{m,1,k})\otimes\cH^{l-p+q}(G_{m,1,r-k})& \text{ if $p+q$ odd}.\end{cases}
\eern

{\rm (iii)} Suppose that $G=SO_N$ and $\max(p,q)\neq 0$. A set of representatives of $W_\fa$-orbits in ${}^0\widehat I$ is 
$$\text{\{$\chi_k$, $0\leq k\leq r/2$\}} \text{ if $\on{min}(p,q)=0$,}\ \{\chi_k,\ 0\leq k\leq r\}\text{ if $p,q\geq 1$}.
$$
 We have 
$$W_{\fa,\chi_k}^0\cong G_{m,1,k}\times G_{m,1,r-k},$$ and 
\begin{align*}&W_{\fa,\chi_{r/2}}/ W_{\fa,\chi_{r/2}}^0\cong\bZ/2\bZ&&\text{ if $\min(p,q)=0$ and $k=r/2$}\\
 &W_{\fa,\chi_k}=W_{\fa,\chi_k}^0&&\text{otherwise}.
\end{align*}
Moreover, we have that
\bern
&&\cH_{W_{\fa,\chi_k}^0}\cong\begin{cases}\cH^{l-p+q}(G_{m,1,k})\otimes \cH^{l+p+q}(G_{m,1,r-k})&\text{ if $p+q$ even}\\\cH^{l+p+q}(G_{m,1,k})\otimes\cH^{l-p+q}(G_{m,1,r-k})& \text{ if $p+q$ odd}.\end{cases}
\eern

\end{proposition}

Let $\on{Irr}\cH_{W_{\La,\chi}^0}$ denote the set of irreducible representations of the Hecke algebra $\cH_{W_{\La,\chi}^0}$. For $\rho\in\on{Irr}\cH_{W_{\La,\chi}^0}$, let $$V_{\rho,\chi}= \bC[\widetilde B_{W_\La}]\otimes_{\bC[\widetilde B_{W_\La}^{\chi,0}]}\otimes(\bC_\chi\otimes\rho).$$ Then $V_{\rho,\chi}$ is an irreducible representation of $\widetilde B_{W_\La}$ when $W_{\La,\chi}^0= W_{\La,\chi}$. We write $V_{\rho,\chi}^\delta$ for its non-isomorphic irreducible summands when $W_{\La,\chi}^0\neq W_{\La,\chi}$. 
\begin{theorem}
\label{main theorem}
{\rm (i)} Suppose that $\theta$ is the order $m=2d$ outer automorphism $\theta^{p,q}$ of $SL_N$, $d=2l+1$, $N=dr+p(p+1)+q^2$. The set of cuspidal character sheaves is
\bern
&&\on{Char}^{\on{cusp}}_{K}(\Lg_1)=\left\{\IC(\Lg_1^\rrs,\cL_\pi)\mid \pi\in{}_A\Theta^{p,q}_d\right\},
\eern
\begin{align*}
{}_A\Theta^{p,0}_d&=&&\{V_{\rho_1\otimes\rho_2,\chi_k}\mid\rho_1\in\on{Irr}\cH^{l+p+1}(G_{d,1,k}),\rho_2\in\on{Irr}\cH^{l+p+1}(G_{d,1,r-k}),\,k\in[0,\tfrac{r}{2}],\,\rho_1\not\cong\rho_2 \}\\
&\cup&&\{V_{\rho\otimes\rho,\chi_{\frac{r}{2}}}^\delta\mid\rho\in\on{Irr}\cH^{l+p+1}(G_{d,1,\frac{r}{2}}),\,\delta={\rm I,II}\}\hspace{.5in}\text{ where $V_{\rho_1\otimes\rho_2,\chi_{\frac{r}{2}}}\cong V_{\rho_2\otimes\rho_1,\chi_{\frac{r}{2}}}$,}\\
{}_A\Theta^{p,q}_d&=&&\{V_{\rho_1\otimes\rho_2,\chi_k}\mid\rho_1\in\on{Irr}\cH^{l+p+q+1}(G_{d,1,k}),\rho_2\in\on{Irr}\cH^{l+p-q+1}(G_{d,1,r-k}),\,k\in[0,r]\}\\&&&\hspace{2in}\text{ if $p+q$ is even and $q\geq 1$},\\
{}_A\Theta^{p,q}_d&=&&\{V_{\rho_1\otimes\rho_2,\chi_k}\mid\rho_1\in\on{Irr}\cH^{l+p-q+1}(G_{d,1,k}),\rho_2\in\on{Irr}\cH^{l+p+q+1}(G_{d,1,r-k}),\,k\in[0,r]\}\\&&&\hspace{2in}\text{ if $p+q$ is odd and $q\geq 1$}.
\end{align*}
{\rm{(ii)}} Suppose that $\theta$ is the order $m=2l$ automorphism $\theta^{p,q}$ of $SO_{N}$, $N=mr+p^2+q^2$, $p+q\leq l$. Suppose that $\max(p,q)=a>0$ and $\min(p,q)=b$. The set of cuspidal character sheaves is
\bern
&&\on{Char}_{K}^{\on{cusp}}(\Lg_1) =\{\IC(\Lg_1^\rrs,\cL_\pi)\mid \pi\in{}_{BD}^0\Theta^{p,q}_m\},
\eern
\begin{align*}
{}_{BD}^0\Theta^{p,q}_m&=&&\{V_{\rho_1\otimes\rho_2,\chi_k}\mid\rho_1\in\on{Irr}\cH^{l+a}(G_{m,1,k}),\rho_2\in\on{Irr}\cH^{l+a}(G_{m,1,r-k}),\,k\in[0,\tfrac{r}{2}],\,\rho_1\not\cong\rho_2\}\\
&
\cup&&\{V_{\rho\otimes\rho,\chi_{\frac{r}{2}}}^\delta\mid\rho\in\on{Irr}\cH^{l+a}(G_{m,1,r/2}),\,\delta={\rm I,II}\},\text{ where $V_{\rho_1\otimes\rho_2,\chi_{\frac{r}{2}}}\cong V_{\rho_2\otimes\rho_1,\chi_{\frac{r}{2}}}$}\\
&&&\hspace{3.5in}\text{ if $\min(p,q)=0$},\\
{}_{BD}^0\Theta^{p,q}_m&=&&\left\{V_{\rho_1\otimes\rho_2,\chi_k}\mid\rho_1\in\on{Irr}\cH^{l-b+a}(G_{m,1,k}),\rho_2\in\on{Irr}\cH^{l+b+a}(G_{m,1,r-k}),\,k\in[0,r]\right\}\\
&&&\hspace{3in}\text{ if $p,q\geq 1$ and $p+q$ is even},\\
{}_{BD}^0\Theta^{p,q}_m&=&&\left\{V_{\rho_1\otimes\rho_2,\chi_k}\mid\rho_1\in\on{Irr}\cH^{l+b+a}(G_{m,1,k}),\rho_2\in\on{Irr}\cH^{l-b+a}(G_{m,1,r-k}),\,k\in[0,r]\right\}\\
&&&\hspace{3in}\text{ if $p,q\geq 1$ and $p+q$ is odd}.
\end{align*}
{\rm (iii)} Suppose that $\theta$ is the order $m=2l$ stable automorphism $\theta^{p,q}$ of $Sp_{2n}$, $2n=mr+p(p+1)+{q(q+1)}$, $p+q\leq l-1$. The set of cuspidal character sheaves is
\bern
&&\on{Char}_{K}^{\on{cusp}}(\Lg_1) =\{\IC(\Lg_1^\rrs,\cL_\pi)\mid \pi\in{}_C\Theta^{p,q}_m\},
\eern
\begin{align*}
{}_C\Theta^{p,q}_m&=&&\left\{V_{\rho_1\otimes\rho_2,\chi_k}\mid\rho_1\in\on{Irr}\cH^{l-p+q}(G_{m,1,k}),\rho_2\in\on{Irr}\cH^{l+p+q+1}(G_{m,1,r-k}),\,k\in[0,r]\right\}\\&&&\hspace{3.5in}\text{ if $p+q$ is even}\\
{}_C\Theta^{p,q}_m&=&&\left\{V_{\rho_1\otimes\rho_2,\chi_k}\mid\rho_1\in\on{Irr}\cH^{l+p+q+1}(G_{m,1,k}),\rho_2\in\on{Irr}\cH^{l-p+q}(G_{m,1,r-k}),\,k\in[0,r]\right\}\\
&&&\hspace{3.5in}\text{ if $p+q$ is odd}.\end{align*}

\end{theorem}

\end{document}